\documentclass[12pt]{amsart}
 
\input{xy} 
\xyoption{all}
\usepackage{graphicx}
\usepackage{color}
\usepackage{verbatim}

\newtheorem{theorem}{Theorem} [section]
\newtheorem{prop}[theorem]{Proposition} 
\newtheorem{lemma}[theorem]{Lemma}
\newtheorem{cor}[theorem]{Corollary}
\newtheorem{conjecture}[theorem]{Conjecture}

\theoremstyle{definition}

\newtheorem{example}{Example}

\theoremstyle{remark}
\newtheorem{remark}[theorem]{Remark}

\numberwithin{equation}{section} 
\numberwithin{example}{subsection}

\newcommand\A{{\mathbb A}}
\newcommand\C{{\mathbb C}}
\newcommand\CC{{\mathbb C}}
\newcommand\Chat { {\hat{\C}} } 

\newcommand\N{{\mathbb N}}

\renewcommand\P{{\mathbb P}}
\newcommand\PP{{\mathbb P}}

\newcommand\RR{{\mathbb R}}
\newcommand\Z{{\mathbb Z}}

\newcommand\Q{{\mathbb Q}}

\newcommand\D{{\mathbb D}}

\newcommand\del{\partial}

\renewcommand\phi{\varphi}

\newcommand\iso{\simeq} 

\newcommand\cM{\mathcal{M}}
\newcommand\cC{\mathcal{C}}

\newcommand{\cP}{\mathcal{P}}

\newcommand\Gal{\operatorname{Gal}}

\newcommand\PSL{\mathrm{PSL}}

\newcommand\supp{\operatorname{supp}}   

\renewcommand\Re {\operatorname{Re}}

\renewcommand\Im {\operatorname{Im}} 

\newcommand\Rat  {\mathrm{Rat}} 
\newcommand\M {\mathrm{M}}

\newcommand\MP {\mathrm{MP}}

\newcommand\Pc{\cP^{cm}}
\newcommand\Per {\mathrm{Per}}

\newcommand\kbar {\overline{k}}
\newcommand\Qbar {\overline{\Q}}
\newcommand\kvbar {\overline{k}_v}
\newcommand\ksep {k^{\rm sep}}
\newcommand\Berk {\operatorname{Berk}} 
\newcommand\ord {\operatorname{ord}} 
\newcommand\hhat {\hat{h}}

\begin{document}

\title{Special curves and postcritically-finite polynomials}

\author{Matthew Baker}
\author{Laura De Marco}
\email{mbaker@math.gatech.edu, demarco@uic.edu}

\date{\today}

\begin{abstract}

We study the postcritically-finite maps within the moduli space of complex polynomial dynamical systems.
We characterize rational curves in the moduli space containing an infinite number of postcritically-finite maps, in terms of critical orbit relations, in two settings: (1) rational curves that are polynomially parameterized and (2) cubic polynomials defined by a given fixed point multiplier.  We offer a conjecture on the general form of algebraic subvarieties in the moduli space of rational maps on $\P^1$ containing a Zariski-dense subset of postcritically-finite maps.  Key ingredients in our proofs are an arithmetic equidistribution theorem, complex-analytic tools for analyzing bifurcation measures, and the theory of polynomial decompositions.
\end{abstract} 


\maketitle

\thispagestyle{empty}

\section{Introduction}

For each integer $d\geq 2$, let $\MP^{cm}_d$ denote the moduli space of critically-marked complex polynomials of degree $d$.%
\footnote{The moduli space $\MP^{cm}_d$ is the space of complex polynomials of degree $d$ modulo conjugacy by conformal automorphisms of $\C$.  It is a finite quotient of $\Pc_d \iso \C^{d-1}$, the space of critically-marked, monic and centered polynomials.  Indeed, $\Pc_d$ may be parameterized by tuples $(c_1, \ldots, c_{d-1}, b) \in \C^d$ such that $c_1 + \cdots + c_{d-1} = 0$.  The associated polynomial is given by $f(z) = d \cdot  \int_0^z \prod_i (\zeta - c_i) \, d\zeta + b,$ with critical points at $\{c_1, \ldots, c_{d-1}\}$ and $b = f(0)$.  Conjugating $f$ by $z\mapsto \omega z$ where $\omega^{d-1} =1$ induces an action of the cyclic group $\Z/(d-1)\Z$ on $\Pc_d$ (coordinatewise multiplication by $\omega$), and the moduli space $\MP_d^{cm}$ is the quotient of $\Pc_d$ under this action.}
We are interested in the postcritically-finite (PCF) polynomials within $\MP^{cm}_d$, i.e., those polynomials whose critical points all have a finite forward orbit under iteration.   Such maps play a fundamental role in the theory of polynomial dynamics. The PCF polynomials form a countable and Zariski-dense subset of $\MP^{cm}_d$; see Proposition \ref{Zariski} below.  Our ultimate goal is to characterize algebraic subvarieties of $\MP^{cm}_d$ containing a Zariski-dense subset of PCF maps.  In this paper, we make some concrete steps in this direction, focusing on certain kinds of algebraic curves in $\MP^{cm}_d$.  We also offer Conjecture \ref{conjecture} for the general setting of subvarieties in the space of rational functions. 

\subsection{Statement of main results}

To illustrate the idea, consider the following family of algebraic curves (introduced by Milnor in \cite{Milnor:cubicpoly}) in the space of critically-marked cubic polynomials:
$$\Per_1(\lambda) = \{f \in \MP^{cm}_3: f \mbox{ has a fixed point with multiplier } \lambda\}$$
for each $\lambda \in\C$.  (Recall that the multiplier of a fixed point is simply the derivative of $f$ at the fixed point.)  

\begin{theorem} \label{Per1}
The curve $\Per_1(\lambda)$ contains infinitely many postcritically-finite cubic polynomials if and only if $\lambda=0$.
\end{theorem}

\noindent
The idea of the proof is as follows.  For $\lambda=0$, one critical point is fixed for all $f\in\Per_1(0)$, so there is exactly one ``active" critical point along each irreducible component of $\Per_1(0)$.  By a classical complex dynamics argument, the active critical point must have finite forward orbit for a dense set of parameters in the bifurcation locus, so there are infinitely many PCF polynomials $f\in\Per_1(0)$.  For the converse direction, assume there are infinitely many postcritically-finite maps in $\Per_1(\lambda)$.  Then $\lambda\in\Qbar$, and we apply an arithmetic equidistribution theorem (Theorem \ref{thm:arithmetic_equidistribution}) to conclude that these PCF maps are equidistributed with respect to the bifurcation measure of each bifurcating critical point.  In particular, if $\lambda\not=0$, then the two critical points define the same bifurcation measure along $\Per_1(\lambda)$.  But the two critical points are dynamically independent and must define distinct bifurcation measures, so we conclude that $\lambda=0$.  See \S\ref{cubics} for details.

In general, we expect that an algebraic subvariety $V$ in $\MP^{cm}_d$ contains a Zariski dense subset of PCF maps if and only if $V$ is cut out by critical orbit relations.  Unfortunately, pinning down a precise notion of ``critical orbit relation" is a bit delicate, as we need to take into account the presence of nontrivial symmetries.   In the next result, we provide a precise formulation for polynomially-parameterized curves in the space $\Pc_d\iso \C^{d-1}$, a branched cover of $\MP^{cm}_d$, consisting of all monic and centered polynomials with marked critical points.  We emphasize the equivalence of statements (1) and (4) in Theorem~\ref{main} below.  

In order to state the result, we first need the following definitions.  A {\em holomorphic family} of polynomials is a holomorphic map $V \to \Pc_d$, $t\mapsto f_t$, from a complex manifold $V$ to the space of polynomials.  A {\em marked point} along $V$ is a meromorphic function $a: V \to \P^1$.  The marked point $a$ is said to be {\em passive} if the sequence of functions $\{t \mapsto f_t^n(a(t)): n\geq 1\}$ forms a normal family on $V$; otherwise it is {\em active}. 
In the algebraic setting, where $V$ is quasiprojective and the family $f_t$ is algebraic (and nonconstant), a marked {\em critical} point of $f_t$ is passive if and only if it is persistently preperiodic \cite[Lemma 2.1]{McMullen:families},  \cite[Theorem 2.5]{Dujardin:Favre:critical}.   

\begin{theorem}  \label{main}
Fix $d\geq 2$.  Let  
	$$f_t =  (c_1(t), \ldots, c_{d-1}(t), b(t)) \in \Pc_d$$
be a holomorphic family of polynomials with marked critical points, defined for $t\in\C$, where each coordinate function lies in $\C[t]$.  The following are equivalent:
\begin{enumerate}
\item		$f_t$ is postcritically finite for infinitely many parameters $t$;
\item		for every pair of active critical points $c_i$ and $c_j$, the normalized bifurcation measures are equal;
\item		the connectedness locus for $\{f_t\}$ is equal to 
			$$\M_i =  \left\{t:  \sup_n |f^n_t(c_i(t))| < \infty\right\}$$
		for any choice of active critical point $c_i$;
\item		for every pair of active critical points $c_i$ and $c_j$, there exist a polynomial $h_t(z) \in \C[t,z]$ and integers $k>0$, $n,m\geq 0$, 
		such that 
			$$h_t \circ f^k_t = f^k_t \circ h_t  \qquad \mbox{and}  \qquad f_t^n(c_j(t)) = h_t (f_t^m(c_i(t)))$$
		for all $t$.
\end{enumerate}
\end{theorem}

\noindent
In plain English, the equivalence of (1) and (4) means that there is a Zariski-dense set of parameters $t\in\C$ for which 
$f_t$ is PCF {\em if and only if} there is exactly one active critical orbit, up to symmetries (the $h$ term).  In particular, the critical point $c_i$ has finite orbit for $f_t$ if and only if $c_j$ has finite orbit for $f_t$.  If $\deg_z h = 1$, then $h_t$ must be a symmetry of the Julia set of $f_t$; these were classified in \cite{Beardon:symmetries}.  If $\deg_z h > 1$, then $h_t$ must share an iterate with $f_t$ for all $t$ \cite{Ritt:permutable}; it follows that condition (4) is symmetric in $i$ and $j$.    In \S\ref{examples}, we provide examples of polynomial families $f_t$ satisfiying the conditions of Theorem \ref{main}, and we illustrate how we can use Theorem \ref{main} to conclude that there are only finitely many postcritically-finite maps in certain explicit families.

Theorem \ref{main} is a special case of the following result which concerns marked (but not necessarily critical) points which are simultaneously preperiodic.    

\begin{theorem}  \label{preperiodic}
Let $f_t$ be a 1-parameter family of polynomials of degree $d\geq 2$, parameterized as 
	$$f_t(z) = z^d + b_2(t)z^{d-2} + \cdots + b_d(t)$$
with $b_j(t) \in \C[t]$ for each $j$.  Let  $a_1(t), a_2(t) \in \C[t]$ be a pair of active marked points, and define
	$$S_i := \{t\in\C:  a_i(t) \mbox{ is preperiodic for } f_t\}.$$  
The following are equivalent:
\begin{enumerate}
\item		$|S_1\cap S_2| = \infty$;
\item		$S_1 = S_2$; and
\item		there exist a polynomial $h\in \C[t,z]$ and integers $k>0$, $n,m\geq 0$ such that 	
		$$h_t \circ f^k_t = f^k_t \circ h_t  \qquad \mbox{and}  \qquad f_t^n(a_1(t)) = h_t f_t^m(a_2(t))$$
		for all $t$.  
\end{enumerate}	
\end{theorem}

\noindent
Theorem \ref{preperiodic} is an extension of the results \cite[Theorem 1.1]{BD:preperiodic} and \cite[Theorem 2.3]{Ghioca:Hsia:Tucker}, where stronger hypotheses guaranteed that the symmetries $\{h_t\}$ must be trivial.  The article \cite{GHT:2012} is closely related, showing that (1) $\iff$ (2) for certain families of rational functions.  As in the case of marked critical points, the activity of $a_i(t)$ in Theorem \ref{preperiodic} is equivalent to the statement that $a_i$ is not persistently preperiodic (Proposition \ref{active}).  

The idea behind our proof of Theorem \ref{preperiodic} is as follows.  If we assume condition (3), then (2) follows immediately and (1) follows from Montel's theorem, showing that an active point must have finite orbit at infinitely many parameters $t$.  For the implication (1) $\implies$ (3), we begin by applying an arithmetic equidistribution theorem (Theorem \ref{thm:arithmetic_equidistribution}) that implies an ``almost (2)" statement:  $S_1$ and $S_2$ can differ by at most finitely many elements.  This step, which uses Berkovich analytic spaces in a crucial way, appeared in \cite{Ghioca:Hsia:Tucker} and we refer there for details.  

To complete the proof that (1) implies (3), we use classical techniques from complex analysis to, first, deduce an analytic relation between the orbits of $a_1$ and $a_2$ and, then, promote this to an invariant algebraic relation.  Finally, via recent results of Medvedev-Scanlon \cite{Medvedev:Scanlon}, employing methods of Ritt \cite{Ritt:decompositions} to classify invariant subvarieties for a certain class of polynomial dynamical systems, we may simplify the form of our algebraic relation to the statement of condition (3).  

Theorem \ref{Per1} is not a special case of Theorems \ref{main} and \ref{preperiodic}, because the rational curves $\Per_1(\lambda)$ in $\Pc_3$ are not parameterized by polynomials for $\lambda\not=0$.    

\subsection{Motivation from results in arithmetic geometry}
\label{context}
In arithmetic geometry, there are numerous results which fit into the following paradigm.  One is given a complex
algebraic variety $X$ and a countable Zariski dense collection of ``special'' algebraic points on $X$.  The question
then arises which algebraic subvarieties of $X$ can contain a Zariski dense set of special points.  Usually one knows
a family of ``special subvarieties'' of $X$ which do contain a Zariski dense set of special points, and the problem is
to determine whether an arbitrary subvariety of $X$ containing a Zariski dense set of special points must itself be special.

The canonical example of this paradigm is the ``Manin-Mumford conjecture'', first established by Raynaud \cite{Raynaud:1, Raynaud:2}.  If $X$ is an abelian variety then the torsion points of $X$ are countable and 
Zariski dense, and if $Y$ is a torsion subvariety of $X$
(meaning a translate of an abelian subvariety by a torsion point) then $Y$ contains a dense set of torsion points.
Conversely, Raynaud's theorem asserts that if an algebraic subvariety $Y$ of $X$ contains a Zariski dense set of torsion points, then $Y$ must be a torsion subvariety.
An analogous result when $X$ is an algebraic torus (so that torsion points are algebraic points of $X$ whose coordinates are all roots of unity) was proved by Laurent, and extended to semiabelian varieties by Hindry \cite{Laurent}, \cite{Hindry}.

A more recent (and in general still conjectural) illustration of the special point and special subvariety formalism is the ``Andr{\'e}-Oort conjecture'';  see e.g.  \cite{Andre:finitude}, \cite{Pila:AO}.  If $X$ is a Shimura variety then the CM points form a countable dense set of algebraic points on $X$, and likewise for any Shimura subvariety $Y$ of $X$.   The Andr{\'e}-Oort conjecture asserts conversely that an algebraic subvariety containing a dense set of CM points must be special, i.e., a Shimura subvariety.  
A concrete special case of this conjecture, proved by Andr{\'e}, is that an irreducible algebraic curve $Y$ in $X = \C^2$ containing a Zariski dense set of points whose coordinates are both $j$-invariants of CM elliptic curves must be either
horizontal, vertical, or a modular curve $X_0(N)$. 

Ghioca, Tucker, and Zhang have put forth some conjectural dynamical analogs of the Manin-Mumford conjecture 
\cite{Ghioca:Tucker:Zhang}. The main results and conjectures in the present paper can be thought of as dynamical analogs of
the Andr{\'e}-Oort conjecture.  The Shimura varieties, which for our purposes can be thought of as moduli spaces for abelian varieties with certain additional structure, get replaced by moduli spaces for polynomial dynamical systems, and CM points get replaced by PCF maps.  As in some approaches to the Manin-Mumford and Andr{\'e}-Oort conjectures, equidistribution theorems for Galois orbits of special points play a crucial role in our approach to the dynamical version of these problems.

\subsection{Examples}  \label{examples}
We now provide examples to illustrate Theorem \ref{main}.  The first few are basic examples of families satisfying the conditions of Theorem \ref{main}.  We include examples where the symmetries $h_t$ are necessarily nontrivial.  We conclude with two examples illustrating how Theorem \ref{main} might be used to show that there are only finitely many postcritically-finite maps in a given family.  

\begin{example} \label{deg2} (Infinitely many postcritically-finite maps in degree 2)  In degree 2, there is a unique critical point, so the space $\MP_2^{cm} \iso \Pc_2$ is itself of dimension 1.  The polynomial $f_t(z) = z^2+t$ is postcritically finite if and only if $t$ satisfies the polynomial equation
	$$f^n_t(0) = f^m_t(0)$$
for some $n>m$.  There are infinitely many such $t$; in fact, a simple argument involving Montel's Theorem shows that they accumulate everywhere in the boundary of the Mandelbrot set.  
\end{example}

\begin{example}  (Maps with an automorphism)
Let $f_t(z) = z^3 - 3t^2z$, so $c_1(t)=t$, $c_2(t) = -t$.  The orbits of $c_1$ and $c_2$ are generally disjoint, though they are symmetric by $h_t(z) = -z$.  That is, we have $h_t\circ f_t = f_t \circ h_t$ and
	$$f_t^n(c_1(t)) = h_t (f_t^n(c_2(t)))$$
for all $t$ and any choice of $n\geq 0$.  There are infinitely many postcritically-finite maps in this family.
\end{example}

\begin{example} (Symmetry of the Julia set)
Let $f_t(z) = z^2(z^3-t^3)$.  The Julia set of $f_t$ has a symmetry of order 3, but $f_t$ has no nontrivial automorphisms for $t\not=0$.  Set $\beta = (2/5)^{1/3}$ and choose $\zeta\not=1$ so that $\zeta^3=1$.  Then $f_t$ has a fixed critical point at $c_1(t) = 0$ for all $t$, and the other critical points are $c_2(t) = \beta t$, $c_3(t) = \zeta\beta t$, $c_4(t) = \zeta^2\beta t$.  Then $f_t(\zeta z) = \zeta^2 f_t(z)$ for all $t$, so $h(z) = \zeta z$ commutes with the second iterate $f_t^2$ and
	$$f_t^2(c_3(t)) = \zeta f^2_t(c_2(t)), \quad f_t^2(c_4(t)) = \zeta f_t^2(c_3(t)), \quad \mbox{ and } f_t^2(c_2(t)) = \zeta f_t^2(c_4(t))$$
for all $t$.  There are infinitely many postcritically-finite maps in this family.
\end{example}

\begin{example} (Symmetry $h$ of degree $>1$)
Let $g_t(z) = z^2-t^2$ and $f_t(z) = g_t^2(z) = (z^2-t^2)^2 - t^2$ of degree 4, with $c_1(t) = 0$, $c_2(t) = t$, $c_3(t) = -t$.  None of the critical points are persistently periodic, and there are infinitely many postcritically-finite parameters for the family $f_t$ (being just the second iterate of the quadratic family).  The critical points $c_2$ and $c_3$ land on $c_1$ after one iterate of $g_t$, but their orbits under $f_t$ are disjoint from the orbit of $c_1$ for all $t\not=0$; however, if we set $h(t,z) = g_t(z)$, then $f_t\circ h_t = h_t \circ f_t$ for all $t$, with 
	$$f_t(c_2(t)) = f_t(c_3(t)) \quad \mbox{and} \quad  h_t (c_i(t)) = c_1(t)$$ 
for all $t$ and $i=2,3$.
\end{example}

\begin{figure} [h]
\includegraphics[width=2.9in]{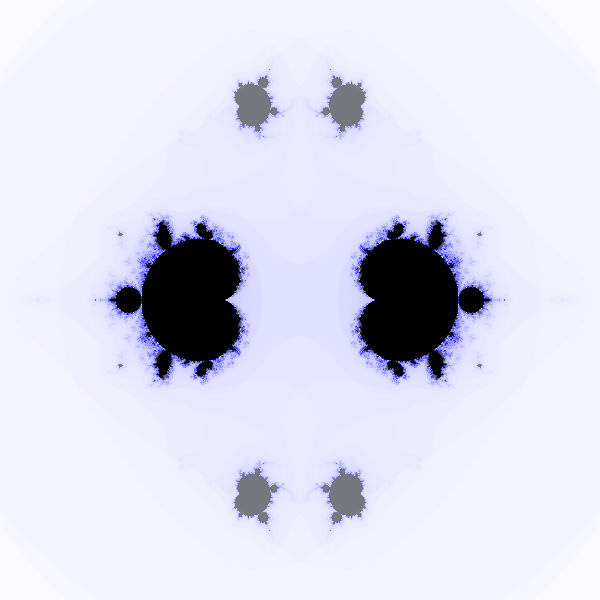}
\includegraphics[width=2.9in]{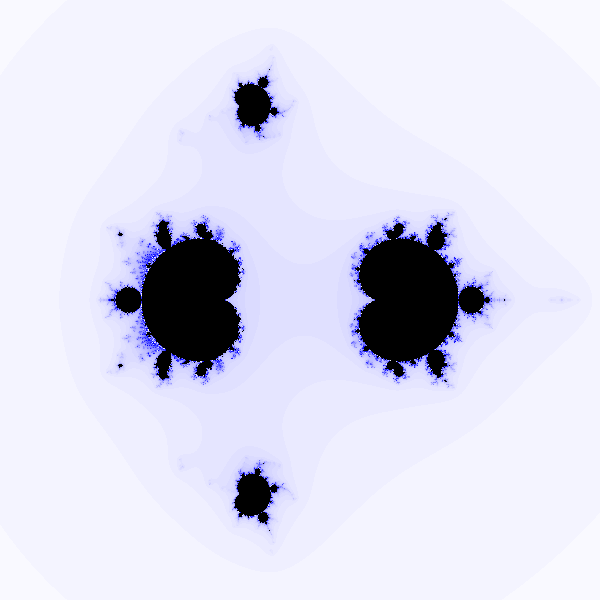} 
\caption{\small Left: the connectedness locus for $f_t(z) = z^3 - 3 t^2 z + 0.56$, of Example \ref{example for figure}, is shown in black in the region $\{|\Re t|, |\Im t| \leq 1.2\}$; gray indicates that only one critical point remains bounded under iteration.  Right:  the boundedness locus $\M_1$ for the critical point $c_1(t) = t$ is shown in black.  The boundedness locus $\M_2$ for $c_2(t) = -t$ is the image of $\M_1$ under $t\mapsto -t$.  The support of the bifurcation measure $\mu_i$ is the boundary of $\M_i$, $i = 1,2$.  }
\label{b=0.56}
\end{figure}

\begin{example} \label{example for figure} (Finitely many PCF polynomials)  
In the example of Figure \ref{b=0.56}, the boundaries of the sets $\M_1$ and $\M_2$ appear to have a great deal of overlap.  Recall that the parameters where critical point $c_i$ has finite forward orbit are dense in the boundary of $\M_i$ (or see Lemma \ref{normal}).  However, there are only finitely many PCF maps, where both critical points have finite forward orbit, by condition (3) of Theorem \ref{main}.  Indeed, there are obvious gray regions in the picture at left, where one critical point remains bounded while the other escapes to infinity.
\end{example}

\begin{example} (Finitely many PCF polynomials)
In the family $f_t(z) = z^3 - 3t^2z + i$, we can employ condition (4) of Theorem \ref{main} to show that (1) fails.  Specifically, if (4) were to hold, the critical point at $t$ would be preperiodic if and only if the critical point at $-t$ is preperiodic.  So it suffices to find a single parameter $t_0$ at which one critical point is preperiodic while the other has infinite forward orbit.  For the parameter $t_0=i$, the critical point at $-i$ is fixed while the critical point at $i$ lies in the basin of infinity.
\end{example}

\subsection{A conjecture for postcritically-finite rational maps}
Let $\{f_t: t\in V\}$ be an $N$-dimensional algebraic family of critically-marked rational maps of degree $d\geq 2$. In other words, $V$ is a quasi-projective complex algebraic variety, the map $t\mapsto f_t$ defines a regular map $V \to \Rat_d \subset \P^{2d+1}_\C$ to the space of rational functions on $\P^1$ of degree $d$, and the image of $V$ in the {\em moduli space} $\M_d$ has dimension $N$, where $\M_d$ is the quotient of $\Rat_d$ by the conjugation action of $\PSL_2\C$.  Furthermore, the critical points of $f_t$ are the images of regular maps 
	$$c_i: V \to \P^1$$
for $i = 1, \ldots, 2d-2$.   A family as above defines a rational function ${\bf f}: \P^1_k \to \P^1_k$ of degree $d$ (where $k=\C(V)$ is the function field of $V$) with critical points $c_i \in \P^1(k)$, $i = 1, \ldots, 2d-2$.

If PCF maps play the role of the ``special points" in the space of rational maps, then the following conjecture provides a characterization of the ``special subvarieties" in the space of critically-marked rational maps $\Rat_d^{cm}$.  An $n$-tuple of marked critical points $(c_{i_1}, \ldots, c_{i_n})$ is said to have {\em dynamically dependent orbits} if there exists a nonzero algebraic relation $\{{\bf P} = 0\}\subset (\P^1_k)^n$, which is invariant under the map $({\bf f}, \ldots,  {\bf f})$, such that  
		$${\bf P}(c_{i_1}, \ldots, c_{i_n}) = 0.$$
Otherwise, we say that the $n$ critical points are {\em dynamically independent} on $V$.  Invariance of $X$ under a map $F$ means that $F(X)\subset X$.   Note that when $n=1$, this definition states that a single critical point is dynamically independent (from itself) if and only if it has infinite orbit for ${\bf f}$; i.e., if and only if it is active.  Moreover, for any $n$, dynamical independence of $n$ critical points forces all $n$ critical points to be active.  

\begin{conjecture}  \label{conjecture}
Suppose $\{f_t: t\in V\}$ is an $N$-dimensional algebraic family of critically-marked rational maps of degree $d\geq 2$, with $V$ irreducible.  Then $f_t$ is postcritically finite for a Zariski-dense subset of $t\in V$ if and only if there are at most $N$ dynamically independent critical points on $V$.
\end{conjecture}

\noindent
In Theorem \ref{main}, our conclusion (4) is stronger than that of Conjecture \ref{conjecture} because we can appeal to the classification results of Medvedev-Scanlon \cite{Medvedev:Scanlon} to obtain a more precise form for the relation ${\bf P}$.  

One implication of Conjecture \ref{conjecture} (dynamical dependence of any active $(N+1)$-tuple of critical points implies Zariski density of PCF maps) follows easily from an argument mimicking the proof of Proposition \ref{Zariski} and the following observation.  If $N+1$ critical points have dynamically dependent orbits along $V$, and if $N$ of them have finite forward orbits at a given parameter $t\in V$, then the $(N+1)$-st critical point will also have finite orbit at $t$.

We remark that the flexible Latt\`es maps, in the case where $d$ is a square, form a 1-dimensional algebraic family of postcritically-finite maps.  Thus, the use of ``at most $N$" in Conjecture \ref{conjecture} (rather than ``exactly $N$") is necessary.  By Thurston's rigidity theorem, this is the only positive-dimensional family with no active critical points; see \cite[Theorem 2.2]{McMullen:families}, \cite{Douady:Hubbard:Thurston}.

\subsection{Acknowledgements}
We had helpful discussions with many people during the course of this project.  Special thanks go to Curt McMullen, Dragos Ghioca, and Tom Tucker for their insightful comments and suggestions.  We thank Alice Medvedev for carefully explaining her work with Tom Scanlon, and Bjorn Poonen for asking the question which led to this work during the Bellairs Workshop on Number Theory in 2010 (``What can one say about subvarieties of the moduli space of rational maps which contain a Zariski-dense set of PCF points?'').  We thank Xander Faber for organizing the Bellairs workshop and ICERM for hosting us during the Spring 2012 Semester in Complex and Arithmetic Dynamics when the details in this paper were solidified.  We also had helpful discussions with Romain Dujardin, Patrick Ingram, and Joe Silverman.  We thank Suzanne Lynch Boyd and Brian Boyd for help with the Dynamics Explorer Tool which was used to gather experimental evidence and generate all of the images in this paper.  Our research was supported by the National Science Foundation.

\bigskip
\section{Activity and normal families}
\label{activity}

In this section we prove the ``easy" implications in Theorem \ref{preperiodic}; the key ingredient is Montel's theory of normal families.  In Proposition \ref{active}, we show that a marked point is passive if and only if it is persistently preperiodic.  We conclude the section with a proof that the PCF polynomials form a countable and Zariski-dense subset of $\Pc_d$ (Proposition \ref{Zariski}).

\subsection{Activity and bifurcation}
Let $f_t$ be a holomorphic family of polynomials of degree $d\geq 2$, parameterized by $t\in\C$.  Let $a: \C\to \C$ be a holomorphic function.  Let 
\begin{equation} \label{escape}
	G_t(z) = \lim_{n\to\infty} \frac{1}{d^n} \log^+ |f_t^n(z)|
\end{equation}
denote the escape-rate function for $f_t$.  Associated to the marked point $a$ is a {\em bifurcation measure} 
\begin{equation} \label{bif measure}
	\mu_a = \frac{1}{2\pi} \Delta G_t(a(t)), 
\end{equation}
where the Laplacian is with respect to $t$, taken in the sense of distributions.  

The name ``bifurcation measure" comes from the special case where $a(t)$ is a critical point of $f_t$ for all $t$.  In that case, the support of $\mu_a$ coincides with the activity locus of the critical point, the set of parameters where the critical point is ``passing through" the Julia set of $f_t$ \cite{Mane:Sad:Sullivan}.  See \cite{D:current} and \cite{Dujardin:Favre:critical} for background on bifurcation currents.  Similarly for any marked point, the support of the measure can be characterized by a bifurcation in its dynamical properties; see e.g., \cite[Theorem 9.1]{D:lyap}.    

Recall that a point $a(t)$ is passive if the sequence of functions $t\mapsto f_t^n(a(t))$, $n\geq 0$, forms a normal family for $t\in\C$; otherwise the point $a$ is active.  In the special case where $a(t)$ is a critical point of $f_t$, the following proposition was established in \cite[Proposition 10.4]{DM:trees} (and a version for rational functions was proved in \cite[Theorem 2.5]{Dujardin:Favre:critical}).  We give a different proof, appealing to properties of the function field height of $f_t$.

\begin{prop}  \label{active}
Let $f_t$ be a family of polynomials, parameterized polynomially as 
	$$f_t(z) = z^d + b_2(t)z^{d-2} + \cdots + b_d(t)$$
with $b_j(t)\in\C[t]$ for each $j$.  Fix a marked point $a(t) \in \C[t]$.  The following are equivalent:
\begin{enumerate}
\item		$a(t)$ is active;
\item		there exists a parameter $t_0\in\C$ for which the forward orbit of $a(t_0)$ under $f_{t_0}$ is infinite;
\item		$G_t(a(t)) = q \log|t| + O(1)$  as $t\to\infty$, for some positive $q\in \Z[1/d]$; and
\item		the bifurcation measure 
			$$\mu_a = \frac{1}{2\pi} \Delta G_t(a(t)) $$
		is nonzero.  
\end{enumerate}
\end{prop}

\proof
We begin with the most delicate implication, that (2) implies (3).  View ${\bf f} = \{f_t\}$ as a polynomial defined over the function field $k=\C(t)$, so ${\bf f} \in k[z]$ and ${\bf a} = a(t) \in \P^1(k)$.  Assuming condition (2), the point ${\bf a}$ is not preperiodic for ${\bf f}$.  Note that ${\bf f}$ is not isotrivial, as $f_t$ is affine conjugate to $f_s$ for only finitely many values of $s$ (where $f_s(z) = \zeta^{-1}f_t(\zeta z)$ with $\zeta^{d-1}=1$).  We may therefore apply \cite[Theorem B]{Benedetto:polynomial} to conclude that the function-field height of ${\bf a}$ is positive.  That is, 
	$$\hat{h}_{\bf f}(a) = \lim_{n\to\infty} \frac{1}{d^n} \log (\deg_t f_t^n(a(t))) > 0$$
so, in particular, 
	$$\deg_tf_t^n(a(t)) \to \infty$$
as $n\to\infty$; see Remark \ref{ff height} for more information.  Choose $n_0$ so that $m_0 = \deg_t f_t^{n_0}(a(t)) > \max_j \deg_t b_j(t)$.  Then for all $n\geq 0$, 
	$$\deg_t f_t^{n+n_0}(a(t)) = m_0d^n.$$
This shows that (2) implies (3) with 
	$$q = \frac{m_0}{d^{n_0}}.$$
	
Condition (1) clearly implies condition (2) (since the field $\C$ is uncountable).  Condition (3) implies condition (4), because the function $G_t(a(t))$ cannot be harmonic on all of $\C$ if it has nontrivial logarithmic growth.  If $\{t\mapsto f_t^n(a(t))\}$ were normal on $\C$, then there would be a subsequence $f_t^{n_k}(a(t))$ that converges locally uniformly in $\C$ to an entire function or to the constant infinity.  But then the escape rate $G_t(a(t))$ would be everywhere 0 or everywhere infinite.  In particular, the measure $\mu_a$ would be trivial.  So (4) $\implies$ (1) and the circuit of implications is closed.  
\qed

\medskip
\begin{remark}  \label{ff height}
We explain briefly the relation between function-field height and degree growth.  Recall that if $k=\C(t)$ with its standard product formula structure and ${\bf f} \in \CC[t,z]$ has degree $d$ as a polynomial in $k[z]$, the 
canonical height $\hat{h}_{\bf f} : \PP^1(\kbar)\to\RR_{\geq 0}$ is defined for $a \in \C[t]$ by
\[
\hat{h}_{\bf f}(a) = \lim_{n \to \infty} \frac{1}{d^n} \sum_{v \in \cM_k} \log^+ |{\bf f}^n(a)|_v.
\]
We can identify $\cM_k$ with $\CC \cup \{ \infty \}$.  For an absolute value $v$ 
corresponding to a point $z \in \CC$, we have $\log^+ |{\bf f}^n(a)|_v = 0$ since 
$\log |{\bf f}^n(a)|_v = -\ord_z({\bf f}^n(a)) \leq 0$.
For $v$ corresponding to the point at infinity,
we have $\log^+|{\bf f}^n(a)|_v = \log |{\bf f}^n(a)|_\infty = \deg({\bf f}^n(a)) \geq 0$.
Thus 
\[
\hat{h}_{\bf f}(a) = \lim_{n \to \infty} \frac{1}{d^n} \deg({\bf f}^n(a)).
\]
\end{remark}

\begin{remark}  \label{M_a compact}  
When the conditions of Proposition \ref{active} are satisfied, the measure $\mu_a$ will be compactly supported in the parameter space $\C$.  Indeed, the function $G_t(a(t))$ is necessarily harmonic where it is positive, as it is a locally-uniform limit of harmonic functions.  The set 
	$$\M_a = \{t\in\C:  \sup_n |f_t^n(a(t))| < \infty\} = \{t\in \C:  G_t(a(t)) = 0\}$$
will be compact.  Up to a multiplicative constant (namely, the $q$ of condition (3)), $t\mapsto G_t(a(t))$ defines the Green's function for $\M_a$ with respect to infinity, and $\mu_a$ (up to scale) is the harmonic measure of $\M_a$ with respect to infinity.
\end{remark}

\subsection{Normality and preperiodic points}
Using Montel's theory of normal families, it is straightforward to prove that the conditions of Proposition \ref{active} guarantee infinitely many parameters for which $a(t)$ has finite forward orbit.  For a proof of Montel's theorem, see \cite[\S3]{Milnor:dynamics}.

\begin{lemma}  \label{normal}
Suppose $f: \D\times\C\to \C$ defines a holomorphic family of polynomials of degree $d\geq 2$, parameterized by the unit disk $\D$.   Let $a: \D\to\C$ be an active marked point.  Then there exists a sequence of distinct parameters $t_n\in\D$ for which $a(t_n)$ is preperiodic for $f_{t_n}$ for all $n\in\N$.  In fact, we can choose the parameter $t_n$ so that $a(t_n)$ lands on a repelling cycle of $f_{t_n}$ for each $n$.  
\end{lemma}

\proof
Let $U$ be the largest open set in $\D$ on which $\{t\mapsto f_t^n(a(t))\}_{n\geq 1}$ is normal; it might be empty, and by the definition of active we know that $U\not=\D$.  Choose $t_0\in \D\setminus U$, and let $\{p_1(t_0), p_2(t_0), \ldots, p_r(t_0)\}$ be any repelling cycle for $f_{t_0}$ of period $r >1$.  By the implicit function theorem, the repelling cycle persists for $t$ in a small neighborhood of $t_0$; let $p_i(t)$ denote the $i$-th point in the corresponding repelling cycle for $f_t$.  Note, in particular, that $p_1(t)\not= p_2(t)$ for all $t$ near $t_0$.   The failure of normality on $\D$ and Montel's Theorem imply there exist a parameter $t_1\in\D$ and an integer $k>1$ such that 
	$$f_{t_1}^k(a(t_1)) \in \{p_1(t_1), p_2(t_1)\}.$$
That is, the point $a(t_1)$ is preperiodic for $f_{t_1}$ and the cycle it lands on is repelling.  Now we repeat the argument:  choose any repelling cycle for $f_{t_0}$ of period $r_2 > r$ and follow it holomorphically in a small neighborhood of $t_0$.  We obtain a parameter $t_2$ so that $a(t_2)$ lands on a repelling cycle for $f_{t_2}$.  As $f_{t_0}$ has repelling cycles of arbitrarily high period, we may repeat the argument indefinitely. By induction, we obtain a sequence $\{t_1, t_2, t_3\ldots\}$ of parameters where $a(t_n)$ is preperiodic for $f_{t_n}$, and for each $n$, $a(t_n)$ lands on a repelling cycle of period $r_n > r_{n-1}$.  
\qed

\begin{prop}  \label{preperiodic easy}
Let $f_t$ be a 1-parameter family of polynomials as in Theorem \ref{preperiodic}, and suppose that active points $a_1(t), a_2(t)\in\C[t]$ satisfy condition (3) of the theorem.  Then both conditions (1) and (2) are satisfied.
\end{prop}

\proof
Because $h_t$ commutes with the iterate $f_t^k$ for all $t$, condition (3) implies immediately that $a_1$ has finite orbit for $f_t$ if and only if $a_2$ has finite orbit for $f_t$.  Thus, condition (2) holds.  For condition (1), it suffices to show that the orbit of $a_1(t)$ is finite for infinitely many parameters $t$.  This is an immediate consequence of Lemma \ref{normal}.  
\qed

\subsection{Countability and density of PCF maps}  \label{PCF background}  
To conclude this section, we provide a proof that the set of PCF maps forms a countable and Zariski dense subset of 
the moduli space of (critically-marked) polynomials of degree $d$.  A sketched proof of density appears in \cite[Proposition 6.18]{Silverman:moduli}, based on the transversality results of Adam Epstein (as appearing in \cite{Buff:Epstein:PCF}), for the corresponding statement in the space of all rational functions of degree $d$.  We provide a more direct argument for density here, from the equivalence of persistently preperiodic and normality of iterates, as first appeared in \cite[Lemma 2.1]{McMullen:families}.   A similar proof shows that PCF maps are Zariski dense in the moduli space of rational maps. The argument that the set of PCF maps is countable (after excluding the flexible Latt\`es maps) requires Thurston's rigidity theorem in the case of rational maps, while we can appeal to compactness of the connectedness locus for polynomials.

\begin{prop}  \label{Zariski}
The PCF polynomials form a countable, Zariski dense subset of $\MP^{cm}_d$.  The coordinates of each PCF polynomial in $\Pc_d$ lie in $\Qbar$.
\end{prop}

\proof
It is convenient to work in the space $\Pc_d\iso \C^{d-1}$, a branched cover of $\MP_d^{cm}$ of degree $d-1$.  A postcritically-finite polynomial $f \in \Pc_d$ is a solution to $d-1$ equations of the form
	$$f^{n_i}(c_i) = f^{m_i}(c_i)$$
for integers $n_i < m_i$, $i = 1, \ldots, d-1$.  As equations in the coordinates of $\Pc_d$, they are polynomials defined over $\Q$. Each postcritically-finite polynomial has connected Julia set; and the connectedness locus is compact in $\Pc_d$ \cite[Corollary 3.7]{Branner:Hubbard:1}.  Consequently, the PCF maps form a countable union of algebraic sets, each contained in a compact subset of $\Pc_d$.  As any compact affine variety is finite, the collection of PCF maps is countable, and each is defined over $\Qbar$.  

We now show Zariski density.  Let $S$ be any proper algebraic subvariety of $\Pc_d$, and let $\Lambda$ be its complement.  It suffices to show that there exists a PCF polynomial in $\Lambda$.  Consider the critical point $c_1$.  
Either it is preperiodic along the quasiprojective variety $\Lambda$ or it is active; see \cite[Lemma 2.1]{McMullen:families} or \cite[Theorem 2.5]{Dujardin:Favre:critical}.  In either case, by applying Montel's theorem if active (as in Lemma \ref{normal} above), there exists a parameter $\lambda_1\in \Lambda$ where $c_1$ is preperiodic.  Suppose $c_1$ satisfies the equation $f^{n_1}(c_1) = f^{m_1}(c_1)$ at the parameter $\lambda_1$.  Let $\Lambda_1 \subset \Lambda$ be the subvariety defined by this equation.  Then $\Lambda_1$ is a nonempty quasiprojective variety, of codimension $\leq 1$ in $\Pc_d$, and $c_1$ is persistently preperiodic along $\Lambda_1$.  

We continue inductively.  Suppose $\Lambda_k$ is a quasiprojective subvariety in $\Pc_d$ of codimension $\leq k$ on which $c_1, \ldots, c_k$ are persistently preperiodic.  If $c_{k+1}$ is persistently preperiodic along $\Lambda_k$, set $\Lambda_{k+1} = \Lambda_k$.  If not, apply Lemma \ref{normal} to find a parameter $\lambda_{k+1}\in \Lambda_k$ where $c_{k+1}$ is preperiodic, and define $\Lambda_{k+1}\subset \Lambda_k$ by the critical orbit relation satisfied by $c_{k+1}$ at $\lambda_{k+1}$.  Then $\Lambda_{k+1}$ has codimension at most $k+1$ in $\Lambda$, and the first $k+1$ critical points are persistently preperiodic along $\Lambda_{k+1}$.  In particular, $\Lambda_{d-1}$ is a nonempty subset of $\Lambda$ and consists of PCF polynomials.
\qed

\begin{remark}
Another proof of the Zariski density in Proposition \ref{Zariski} follows from the following theorem of Dujardin and Favre:  the closure of the set of postcritically-finite polynomials (in the usual analytic topology) contains the support of the bifurcation measure in $\MP_d^{cm}$ \cite[Corollary 6]{Dujardin:Favre:critical}.  The bifurcation measure $\mu_{\rm bif}$ cannot charge pluripolar sets \cite[Proposition 6.11]{Dujardin:Favre:critical}, and so the PCF maps are Zariski dense.     
\end{remark}

\bigskip
\section{Arithmetic equidistribution}
\label{equidistribution}

In this section we recall a general arithmetic equidistribution theorem which will be used in the sequel.
We state this result in a form which is a hybrid of the terminology from \cite{BRBook} and \cite{FRL:equidistribution}; the proof follows directly from the arguments in either of those works.\footnote{A closely related equidistribution theorem was proved independently by Chambert-Loir \cite{ChambertLoir}.}  
The result is most naturally formulated using Berkovich spaces; see \cite{BRBook} for an overview.

Let $k$ be a {\em product formula field}.  This means that $k$ is equipped with a set $\cM_k$ of pairwise inequivalent nontrivial absolute values, together with
a positive integer $N_v$ for each $v \in \cM_k$, such that
\begin{itemize}
\item for each $\alpha \in k^\times$, we have $|\alpha|_v = 1$ for all but finitely many $v \in \cM_k$; and
\item  every $\alpha \in k^\times$ satisfies the {\em product formula}
\index{product formula}%
\begin{equation*} 
\prod_{v \in \cM_k} |\alpha|_v^{N_v} \ = \ 1 \ .
\end{equation*}  
\end{itemize}
Examples of product formula fields are number fields and function fields of normal projective varieties.

Let $\kbar$ (resp. $\ksep$) denote a fixed algebraic (resp. separable) closure of $k$. 
For $v \in \cM_k$, let $k_v$ be the completion of $k$ at $v$, 
let $\kvbar$ be an algebraic closure of $k_v$, 
and let $\CC_v$ denote the completion of $\kvbar$.  
For each $v \in \cM_k$, we fix an embedding of $\kbar$ in $\CC_v$ extending the canonical embedding of $k$ in $k_v$.
For each $v \in \cM_k$, we let $\PP^1_{\Berk,v}$ denote the Berkovich projective line over $\CC_v$, which is a canonically defined path-connected compact Hausdorff 
space containing $\PP^1(\CC_v)$ as a dense subspace.
If $v$ is Archimedean, then $\CC_v \cong \CC$ and $\PP^1_{\Berk,v} = \PP^1(\CC)$.

For each $v \in \cM_k$ there is a naturally defined distribution-valued Laplacian operator $\Delta$ on $\PP^1_{\Berk,v}$.
For example, the function $\log^+|z|_v$ on $\PP^1(\CC_v)$ extends naturally to a continuous real valued function $\PP^1_{\Berk,v} \backslash \{ \infty \} \to \RR$ and
\[
\Delta \log^+|z|_v = \delta_{\infty} - \lambda_v,
\]
where $\lambda_v$ is the uniform probability measure on the complex unit circle $\{ |z| = 1 \}$ when $v$ is archimedean and 
$\lambda_v$ is a point mass at the Gauss point of $\PP^1_{\Berk,v}$ when $v$ is non-archimedean.

A probability measure $\mu_v$ on $\PP^1_{\Berk,v}$ is said to have {\em continuous potentials} if $\mu_v - \lambda_v = \Delta g$
with $g : \PP^1_{\Berk,v} \to \RR$ continuous.
If $\mu$ has continuous potentials then there is a corresponding {\em  Arakelov-Green function} $g_{\mu} : \PP^1_{\Berk,v} \times \PP^1_{\Berk,v} \to \RR \cup \{ +\infty \}$ which is characterized by the differential equation $\Delta_x g_{\mu}(x,y) = \delta_y - \mu$ and the normalization $\iint g_{\mu}(x,y) d\mu(x) d\mu(y) = 0$.
The function $g_{\mu}$ is finite-valued and continuous outside of 
\[
{\rm Diag_v} := \{ (z,z) \in \P^1(\CC_v) \times \P^1(\CC_v) \} \subseteq \PP^1_{\Berk,v} \times \PP^1_{\Berk,v}.
\]

If $\rho,\rho'$ are measures on $\PP^1_{\Berk,v}$ and $\mu=\mu_v$ is a probability measure with continuous potentials, we define the 
{\em $\mu$-energy} of $\rho$ and $\rho'$ by
\[
( \rho, \rho' )_{\mu} := \frac{1}{2} \iint_{\PP^1_{\Berk,v} \times \PP^1_{\Berk,v} \backslash {\rm Diag}} g_{\mu}(x,y) d\rho(x) d\rho'(y).
\]
One can show that if $\rho$ and $\rho'$ have total mass zero then
$((\rho, \rho')) := ( \rho, \rho' )_{\mu}$ is independent of $\mu$; in this case our definition and notation coincide with those of Favre and Rivera-Letelier \cite{FRL:equidistribution}.

An {\em adelic measure} on $\PP^1$ (with respect to the product formula field $k$) is a collection ${\mathbb \mu} = \{ \mu_v \}_{v \in M_k}$ of probability measures 
on $\PP^1_{\Berk,v}$, one for each $v \in \cM_k$, such that
\begin{itemize}
\item  $\mu_v = \lambda_v$ for all but finitely many $v \in \cM_k$; and
\item  $\mu_v$ has continuous potentials for all $v \in \cM_k$.
\end{itemize}

For a finite subset $S$ of $\PP^1(\ksep)$ and $v \in \cM_k$, we denote by $[S]_v$ the discrete probability measure on $\PP^1_{\Berk,v}$ supported 
equally on all elements of the $\Gal(\ksep / k)$-orbit of $S$.
The {\em canonical height} of $S$ associated to the adelic measure ${\mathbb \mu}$ is defined by
\[
\hhat_{{\mathbb \mu}}(S) := \sum_{v \in \cM_k} N_v \cdot ([S]_v,[S]_v)_{\mu_v}.
\]
(For a justification of the term `canonical height', see for example \cite[Lemma 10.27]{BRBook}.)
This is a Weil height function, in the sense that there is a constant $C$ such that $|h(z) - \hhat_{{\mathbb \mu}}(z)| \leq C$
for all $z \in \ksep$, where $h$ is the standard logarithmic height on $\PP^1$.

\begin{theorem} \cite{BRBook, FRL:equidistribution}
\label{thm:arithmetic_equidistribution}
Let $\hhat_{{\mathbb \mu}}$ be the canonical height associated to an adelic measure $\mu$.
Let $\{S_n\}_{n\geq 0}$ be a sequence of finite subsets of $\PP^1(\ksep)$ for which 
	$$\# (\Gal(\ksep / k) \cdot S_n) \to \infty \qquad \mbox{and} \qquad \hhat_{{\mathbb \mu}}(S_n) \to 0$$ 
as $n \to \infty$.
Then $[S_n]_v$ converges weakly to $\mu_v$ on $\PP^1_{\Berk,v}$ as $n \to \infty$ for all $v \in \cM_k$.
\end{theorem}

\begin{remark}
When $k$ is a number field, Theorem~\ref{thm:arithmetic_equidistribution} is essentially the same as Theorem 2 of \cite{FRL:equidistribution}.
Special cases of Theorem~\ref{thm:arithmetic_equidistribution}, for arbitrary $k$, are proved in Theorems 7.52 and 10.24 of \cite{BRBook}.
It is straightforward to prove the general case of Theorem~\ref{thm:arithmetic_equidistribution} (for arbitrary $k$) by using 
\cite[Lemma 7.55]{BRBook} in conjunction with the proof of \cite[Theorem 2]{FRL:equidistribution}, as in the proof of \cite[Theorem 7.52]{BRBook}.
\end{remark}

\begin{remark} 
\label{NorthcottRemark}
If $k$ is a number field and $S_n$ is the set of $\Gal(\kbar/k)$-conjugates of $z_n$, 
then $\# S_n \to \infty$ follows automatically from the assumption that $\hhat_{{\mathbb \mu}}(S_n) \to 0$
by Northcott's theorem and the fact that $h_{{\mathbb \mu}}$ is a Weil height.
\end{remark}



In order to apply Theorem~\ref{thm:arithmetic_equidistribution} in practice, one often needs to know how to explicitly compute the Arakelov-Green's functions $g_{\mu_v}(x,y)$ for $v \in \cM_k$.  There is a particularly nice way to do this when each $\mu_v$ is the equilibrium measure of a compact set $E_v \subset \A^1_{\Berk,v}$, which will always be the case for the applications in the present paper.  In order to explain how this works, we introduce some terminology.

Fix a place $v$ of $k$ and suppose that $\mu_v$ is the equilibrium measure for a compact set $E_v \subset \A^1_{\Berk,v}$.
Let $G_v : \A^1_{\Berk,v} \to \RR$ be the Green's function for $E_v$, which by assumption is continuous (i.e., we assume that $E_v$ is a {\em regular} set). 
Let $\gamma_v$ be the Robin constant of $E_v$, so the logarithmic capacity of $E_v$ is $e^{-\gamma_v}$ and 
$G_v(s) = \log|s|_v + \gamma_v + o(1)$ as $s \to \infty$.

Define $H_v : \CC_v^2 \to \RR$ by
\[
H_v(s,t) = \left\{ 
\begin{array}{ll}
G_v(s/t) + \log|t|_v & t \neq 0 \\
\log|s|_v + \gamma_v & t = 0. \\ 
\end{array} \right.
\]
Then $H_v$ is continuous and scales logarithmically, i.e., $H_v(\alpha s, \alpha t) = H_v(s,t) + \log|\alpha|_v$.

The following formula comes from a straightforward calculation which we omit.

\begin{prop}
\label{prop:applying_equi}
The normalized Arakelov-Green function $g_{\mu_v}(x,y)$ with respect to $\mu_v$ is given, for $x,y \in \PP^1(\CC_v)$, by the explicit formula
\begin{equation}
\label{eq:explicitgmu}
g_{\mu_v}(x,y) = -\log|\tilde{x} \wedge \tilde{y}|_v + H_v(\tilde{x}) + H_v(\tilde{y}) - \gamma_v,
\end{equation}
where $\tilde{x},\tilde{y}$ are arbitrary lifts of $x,y$ to $\CC_v^2 \setminus \{ 0 \}$ and $(s_1,t_1)\wedge (s_2,t_2) = s_1 t_2 - s_2 t_1$.
\end{prop}

\begin{remark}
For $v$ archimedean, the fact that $g_{\mu_v}(x,y)$ is normalized implies (and in fact is equivalent to) the statement that 
$e^{-\gamma_v}$ is the homogeneous capacity (in the sense of \cite{D:lyap}) of the set $K = \{(s,t)\in\C^2:  H \leq 0\}.$
This is proved in a slightly more roundabout way in \cite[\S4]{D:lyap}.
\end{remark}

Applying the product formula to (\ref{eq:explicitgmu}), we obtain:

\begin{cor}
\label{cor:applying_equi}
Let ${\mathbb \mu} = \{ \mu_v \}_{v \in \cM_k}$ be an adelic measure such that $\mu_v$ is the equilibrium measure associated to a compact set $E_v \subset \A^1_{\Berk,v}$
for all $v \in \cM_k$.  Assume that the global Robin constant $\gamma := \sum N_v \gamma_v$ is zero.  Let $S \subset k$ be a $\Gal(\ksep / k)$-stable finite set such that $G_v(z)=0$ for every $v \in \cM_k$ and every $z \in S$.  Then $\hhat_{{\mathbb \mu}}(S) = 0$.
\end{cor}

\bigskip
\section{Cubic polynomials and fixed point multipliers}
\label{cubics}

Our goal in this section is to prove Theorem \ref{Per1}.  We begin with the ``easy" implication, as in \S\ref{activity}, and prove in Proposition \ref{Per1 easy} that there are infinitely many PCF maps in $\Per_1(0)$. We then show that for each $\lambda\not=0$, there are only finitely many conformal conjugacy classes of postcritically-finite polynomials with a fixed point of multiplier $\lambda$.  For this implication, we apply the arithmetic equidistribution results described in \S\ref{equidistribution}. 

Though our proof of Theorem \ref{Per1} does not use this, we remark that it suffices to study the curves for $|\lambda|>1$.  Indeed, if $0<|\lambda|\leq 1$, there are no postcritically-finite maps on $\Per_1(\lambda)$; see e.g., \cite[Corollary 14.5]{Milnor:dynamics}.

\subsection{The case of $\lambda=0$}  
\begin{figure} [h]
\includegraphics[width=6.0in]{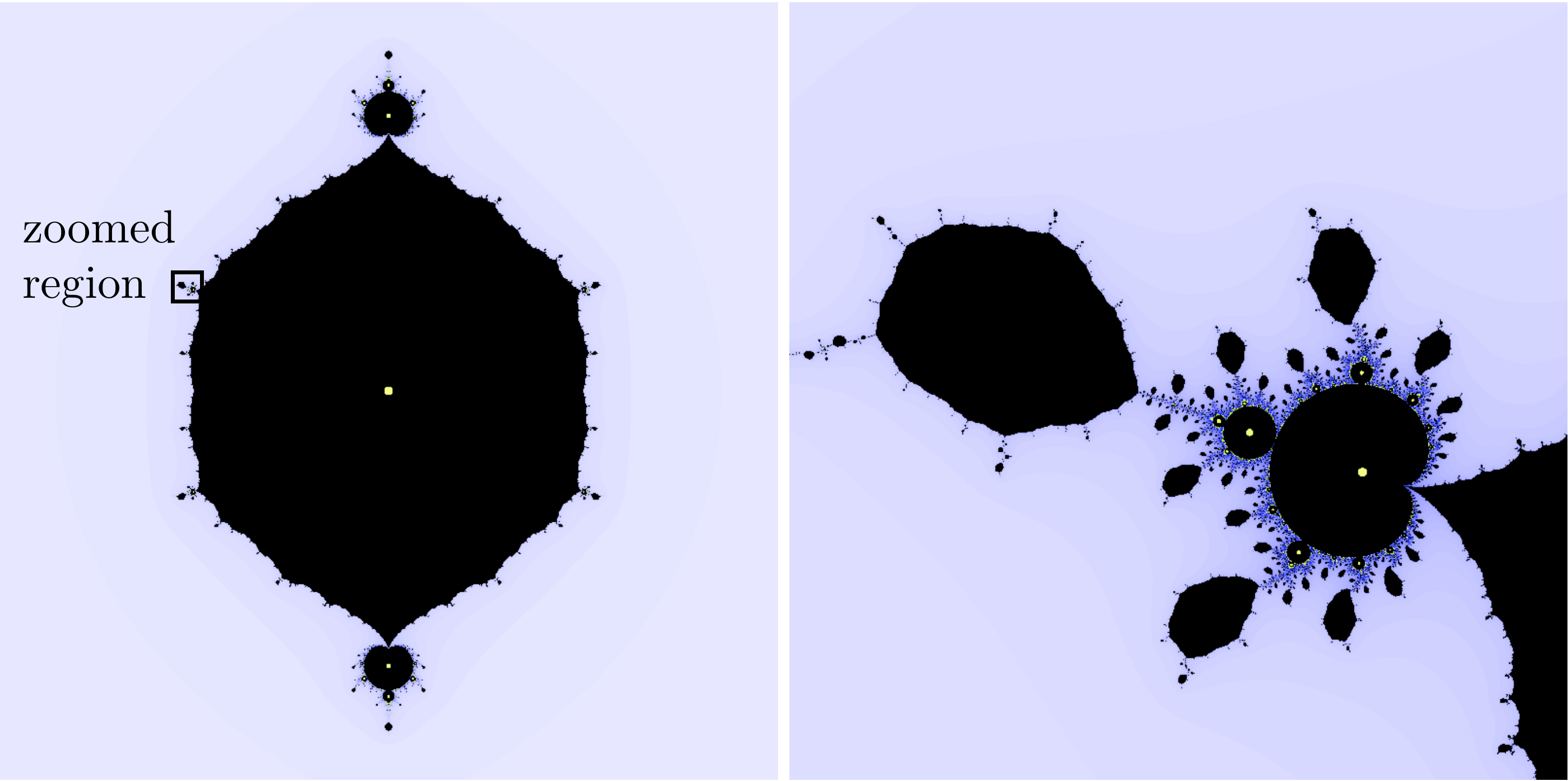} 
\caption{\small At left:  an illustration of the connectedness locus in one component of $\Per_1(0)$ in $\Pc_3$.  The parameters where the active critical point is periodic are marked in yellow.  Proposition \ref{S equidistribution} shows that these parameters are equidistributed with respect to the bifurcation measure.  The curve may be parametrized as $f_t(z) = z^3 - 3t^2 z + t + 2t^3$, where the critical point $c_1(t) = t$ is fixed for all $t$, while $c_2(t) = -t$ is active.   At right:  a zoom of one of the small copies of the Mandelbrot set.}
\label{Per1_0}
\end{figure}

By definition, the curve $\Per_1(0)$ consists of all conjugacy classes of cubic polynomials for which one critical point is fixed.  See Figure \ref{Per1_0}.

\begin{prop}  \label{Per1 easy}
There are infinitely many postcritically-finite cubic polynomials in $\Per_1(0) \subset\MP^{cm}_3$. 
\end{prop}

\proof
It is convenient to work in the space $\Pc_3\iso \C^2$, which is a degree-2 branched cover of $\MP_3^{cm}$.  In $\Pc_3$, the curve $\Per_1(0)$ has two irreducible components, and througout each component, one of the two marked critical points is fixed.  Recalling that the connectedness locus $\cC_3 = \{f\in\Pc_3: J(f) \mbox{ is connected}\}$ is compact \cite[Corollary 3.7]{Branner:Hubbard:1}, we see that both of the critical points cannot be persistently preperiodic along a component of $\Per_1(0)$.  Indeed, one critical point must escape to infinity (and therefore have infinite orbit) for parameters outside the connectedness locus.  Thus, exactly one critical point is active on each component of $\Per_1(0)$.  A polynomial in $\Per_1(0)$ is postcritically-finite if the active critical point has finite forward orbit.  By Lemma \ref{normal}, there are infinitely many postcritically-finite polynomials in $\Per_1(0)$.  
\qed

\subsection{Parameterization of $\Per_1(\lambda)$}  \label{parametrization}
Fix $\lambda \in \C\setminus\{0\}$.  To study the curve $\Per_1(\lambda)$ in the moduli space of cubic polynomials with marked critical points, we shall work with the following parameterization:
	$$f_s(z) = \lambda\, z - \frac{\lambda}{2} \left( s + \frac{1}{s} \right) z^2 + \frac{\lambda}{3} \, z^3$$
for $s\in \C\setminus\{0\}$.  The polynomial $f_s$ has a fixed point at $z=0$ with multiplier $\lambda$ and critical points at $c_+(s) = s$ and $c_-(s)=1/s$.  It is conjugate to the centered polynomial 
	$$P_s(z) = \frac{\lambda}{3} \, z^3 + \left(\frac{\lambda}{2} - \frac{\lambda}{4}\left(s^2 + \frac{1}{s^2}\right) \right)  z +
		\frac{1}{12}  \left(s + \frac{1}{s}\right) \left( 6 - 4\lambda + \lambda s^2 + \frac{\lambda}{s^2} \right)$$
with critical points at $\pm (s^2-1)/(2s)$.  Therefore, the family $f_s$ projects to the curve $\Per_1(\lambda)$ within $\Pc_3$ via 
	$$s\mapsto \left( \sqrt{\frac{\lambda}{3}} \; \frac{s^2-1}{2s}, - \sqrt{\frac{\lambda}{3}} \; \frac{s^2-1}{2s} , - \frac{1}{12} \sqrt{\frac{\lambda}{3}} \left(s + \frac{1}{s}\right) \left( 6 - 4\lambda + \lambda s^2 + \frac{\lambda}{s^2} \right) \right)$$
for either choice of $\sqrt{\lambda}$.   This projection is generically one-to-one.  This curve in $\Pc_3$ then projects to $\Per_1(\lambda)$ in $\MP_3^{cm}$ with degree two, via the identification of $(c_1(s), c_2(s), b(s))$ with $(- c_1(s), -c_2(s), -b(s)) = (c_2(s), c_1(s), -b(s))$.

\subsection{The bifurcation measures}  \label{mu+mu-}
Consider the escape-rate functions
	$$G^+(s) = \lim_{n\to\infty} \frac{1}{3^n} \log^+ |f_s^n(s)|$$
and
	$$G^-(s) = \lim_{n\to\infty} \frac{1}{3^n} \log^+ |f_s^n(1/s)|.$$
An induction argument shows immediately that $f_s^n(s)$ is a polynomial in $s$ for all $n$.  In fact, 
	$$f_s^n(s) = \frac{\lambda}{3} \left( \frac{\lambda}{3} \right)^3 \cdots \left(\frac{\lambda}{3}\right)^{3^{n-2}} 
		\left( \frac{-\lambda}{6} \right)^{3^{n-1}} s^{3^n} +\;  O(s^{3^n-1}) \; \in \; \C[s],$$
so 
\begin{equation} \label{G+}
	G^+(s) = \log|s| + \log|\lambda/3|^{1/6} + \log|\lambda/6|^{1/3} + o(1)
\end{equation}
as $s\to\infty$ and $G^+(s)$ is bounded for $s$ near 0.  By symmetry, $G^-(s) = G^+(1/s)$, so $G^-$ has a logarithmic singularity at $s=0$ and remains bounded as $s\to\infty$.

\begin{lemma}
For each $\lambda\not=0$, both critical points of $f_s$ are active.  
\end{lemma}

\proof
This follows immediately from the nontrivial growth of $G^+$ and $G^-$.
\qed

\bigskip
The bifurcation measures of the critical points $c_+(s) = s$ and $c_-(s)= 1/s$ are defined by 
	$$\mu_+ = \frac{1}{2\pi} \Delta G^+$$
and 
	$$\mu_- = \frac{1}{2\pi} \Delta G^-$$
on $\C\setminus\{0\}$.  From the growth of $G^+$ and $G^-$, we see that $\mu_+$ and $\mu_-$ define probability measures on $\Chat = \C\cup\{\infty\}$.  The support of $\mu_+$ is compactly contained in $\C$, and it does not put positive mass on $s=0$.  Similarly for $\mu_-$.  

\begin{figure} [h]
\includegraphics[width=5.85in]{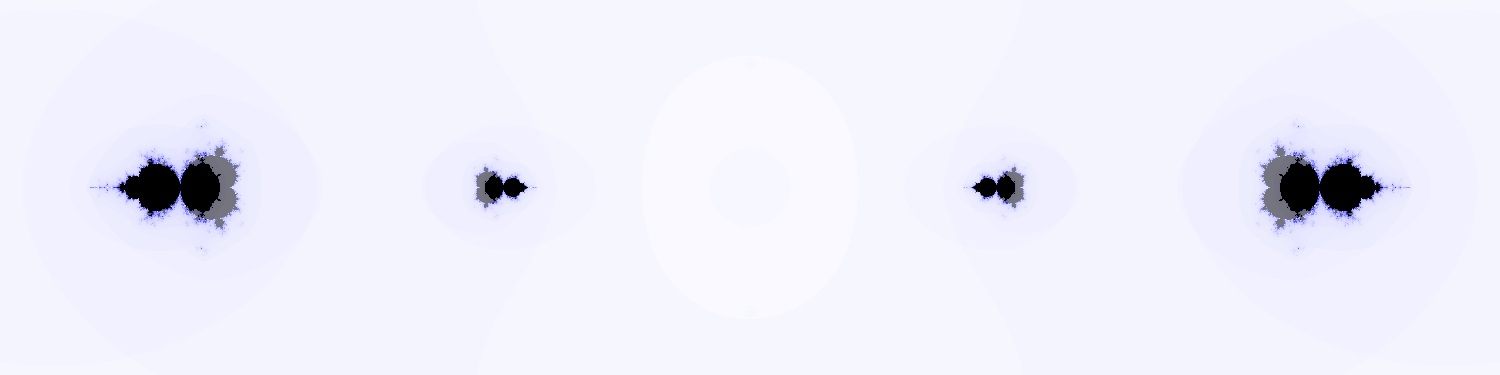}
\includegraphics[width=2.9in]{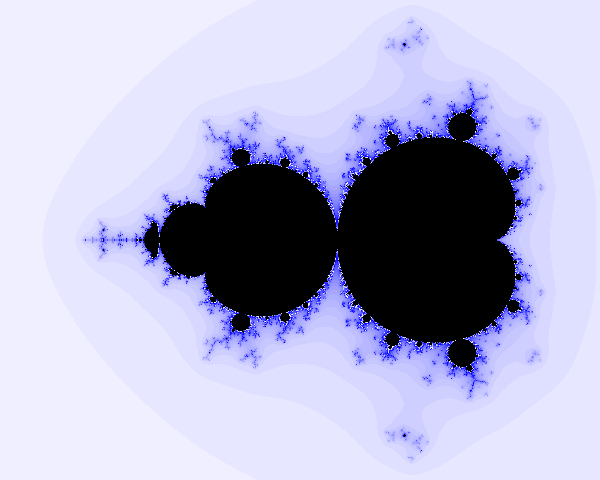} 
\includegraphics[width=2.9in]{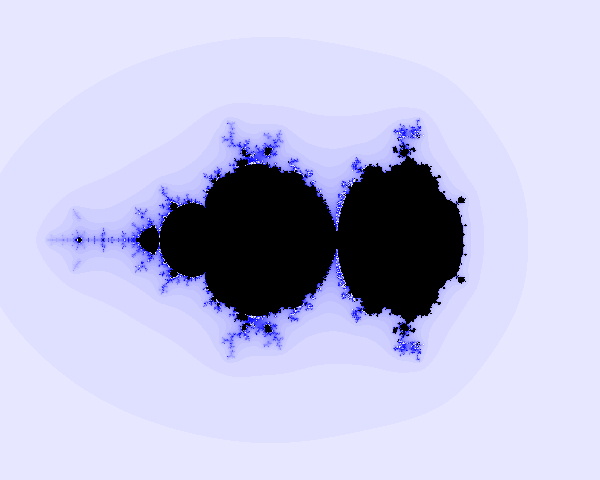} 
\caption{ \small Top:  $\Per_1(6)$ in the parameterization of \S\ref{parametrization}, with $|\Re s| \leq 2$, $|\Im s| \leq 0.5$.  The connectedness locus is shown in black, while gray indicates that only one critical point remains bounded under iteration.    Bottom left:  the support of $\mu_+$ in $\Per_1(6)$, in the region $\{-1.8\leq \Re s \leq -1.3, \;-0.2 \leq \Im s \leq 0.2\}$, is the boundary of the black set, with level sets of $G^+$ shown in shades of blue.  Bottom right: the support of $\mu_-$ and level sets of $G^-$ in the same region.  The polynomial $f_s$ is PCF for $s = -(1+\sqrt{5})/2$, where the two critical points form a cycle of period 2.  As in the example of Figure \ref{b=0.56}, there appears to be a great deal of overlap between the activity locus of $c_+$ and that of $c_-$, though Theorem \ref{Per1} tells us that there are only finitely many PCF maps in $\Per_1(6)$.} 
\label{A6}
\end{figure}

The bifurcation locus for the family $\{f_s\}$ is the set of parameters $s_0$ where the Julia sets $J(f_s)$ fail to vary continuously (in the Hausdorff topology) on any neighborhood of $s_0$.  The bifurcation locus coincides with $(\supp\mu_+) \cup (\supp \mu_-)$; see \cite[Theorem 1.1]{D:current} or \cite[Theorem 3.2]{Dujardin:Favre:critical} for a proof.  See Figure \ref{A6} for an illustration of the bifurcation locus in $\Per_1(6)$.  

We thank Curt McMullen for suggesting the proof of this next lemma: 

\begin{lemma} \label{currents}
The bifurcation measures $\mu_+$ and $\mu_-$ are not equal in $\Per_1(\lambda)$.  
\end{lemma}

\proof
Suppose $\mu_+ = \mu_-$.  Let $B$ denote the bifurcation locus, so $B = \supp \mu_+ = \supp \mu_-$ is compactly contained in $\C\setminus\{0\}$.   The function $G^+ - G^-$ must be harmonic on $\C\setminus\{0\}$, and from the computation of the escape-rate functions above, $G^+-G^-$ grows logarithmically at each end.  Therefore,  $G^+(s) - G^-(s) = C + \log|s|$ for some constant $C$.  Therefore $B = \{G^+ = G^- = 0\} \subset \{G^+ - G^- = 0\}$, so $B$ is a subset of a circle.  But the bifurcation locus $B$ must contain homeomorphic copies of the Mandelbrot set, by the universality of $\del M$ \cite{McMullen:universal}.  This is a contradiction.  
\qed

\subsection{Proof of Theorem \ref{Per1}} \label{proof of Per1}

For $\lambda=0$, Proposition \ref{Per1 easy} states that there are infinitely many postcritically-finite polynomials in $\Per_1(0)$.  

Now suppose $\lambda\not=0$.  
Because all PCF polynomials are defined over $\overline{\mathbb{Q}}$ (Proposition \ref{Zariski}), the existence of a PCF map in $\Per_1(\lambda)$ implies that $\lambda$ is algebraic.  We may therefore assume that the family $f_s$ is defined over a number field $k$.  

We now set up the technical apparatus needed to apply arithmetic equidistribution (Theorem \ref{thm:arithmetic_equidistribution}), to see that parameters where one of the critical points has finite orbit are equidistributed with respect to its bifurcation measure.  We use homogeneous coordinates on both the parameter space and the dynamical space.  For each place $v$ of $k$, let $\C_v$ be the completion of an algebraic closure of the completion of $k$ with respect to $v$, and define
	$$F_{(s,t)}:  \C_v^2 \to \C_v^2$$
by 
$$F_{(s,t)}(z,w) = \left(\lambda\, z w^2 - \frac{\lambda}{2} \left( \frac{s}{t} + \frac{t}{s} \right) z^2w + \frac{\lambda}{3} \, z^3 \, , \; w^3\right)$$
with $(s,t) \in \C_v^*\times\C_v^*$.  Note that $F_{(s,t)} = F_{(t,s)}$ and $f_s(z)$ is the first coordinate of $F_{(s,1)}(z,1)$.  We define 
	$$H_v^+(s,t) = \lim_{n\to\infty} \frac{1}{3^n} \log \|F_{(s,t)}^n(s,t)\|_v$$
and 
	$$H_v^-(s,t) = \lim_{n\to\infty} \frac{1}{3^n} \log \|F_{(s,t)}^n(t,s)\|_v,$$
where $\| (a,b) \|_v = \log \max (|a|_v,|b|_v)$.
Both $H_v^+$ and $H_v^-$ satisfy 
	$$H^\pm_v(\alpha s, \alpha t) = H^\pm_v (s,t) + \log |\alpha|_v$$
for any $\alpha \in \C_v^*$.   

Note that 
\begin{equation}
\label{eq:Gvexpansion}
G_v^+(s) = H_v^+(s,1) = \log|s|_v + \log|\lambda/3|_v^{1/6} + \log|\lambda/6|_v^{1/3} + o(1)
\end{equation}
as $s \to \infty$ by the same calculation as in (\ref{G+}),
and that $G_v^+(s)$ extends continuously to $\A^1_{\Berk,v}$.
Moreover, one sees easily that:
\begin{itemize} 
\item[(G1)] $G_v^+(s)$ is the Green's function relative to $\infty$ for the set 
\[
E_v^+ = \{ z \in \A^1_{\Berk,v} \; : \; G_v^+(z) = 0 \}.
\]
In particular, the Robin constant for $E_v^+$ is $\gamma_v =  \log|\lambda/3|_v^{1/6} + \log|\lambda/6|_v^{1/3}$ by (\ref{eq:Gvexpansion})
and the global Robin constant $\gamma = \sum N_v \gamma_v$ is equal to zero by the product formula.
\item[(G2)] $G_v^+(s)=0$ whenever the polynomial $f_s$ is PCF.
\end{itemize}

Let $\mu_v^+$ be the equilibrium measure for $E_v^+$ (when $v$ is archimedean, this coincides with the probability measure $\mu_+$ introduced in \S\ref{mu+mu-}) 
and let $\mu^+ = \{ \mu_v^+ \}_{v \in \cM_k}$ be the corresponding adelic measure.
(Note that this is indeed an adelic measure, as it is straightforward to verify that $E_v^+$ is the unit disk $\{ z \in \A^1_{\Berk,v} \; : \; |z|_v \leq 1 \}$ in $\A^1_{\Berk,v}$ for all but finitely many places $v$ of $k$.)   Let $\hhat_+$ denote the associated canonical height function.

Let $\{s_n\} \subset \kbar$ be any infinite sequence of parameters such that $c_+(s_n)$ has finite orbit for $f_{s_n}$, for all $n$.  Let $S_n$ denote the set of $\Gal(\kbar/k)$-conjugates of $s_n$.  By (G1) and (G2), the hypotheses of Corollary~\ref{cor:applying_equi} are satisfied, and we conclude that $\hhat_+(S_n)=0$ for all $n$.  Theorem \ref{thm:arithmetic_equidistribution} (combined with Remark~\ref{NorthcottRemark}) applies to show:

\begin{prop}  \label{c_+ equidistribution}  Fix $\lambda\in\Qbar\setminus\{0\}$.  
For any infinite sequence of parameters $s_n\in\kbar$ for which $c_+$ is preperiodic, the discrete probability measures $[S_n]_v$ converge weakly to $\mu_v^+$ in $\A^1_{\Berk,v}$, for all places $v$ of $k$.  In particular, the preperiodic parameters are equidistributed with respect to the bifurcation measure $\mu_+$ on $\C$.
\end{prop}


The same considerations apply to $G_v^-$ and $\mu_v^-$, so Proposition \ref{c_+ equidistribution} holds also for $c_-$ and the adelic measure $\mu^-$.  

Finally, suppose that $\{f_{s_n}\}$ is an infinite sequence of PCF maps in $\Per_1(\lambda)$, so that both critical points are preperiodic for all $n$.  The Galois orbits of $s_n$ must be equidistributed with respect to both $\mu_v^+$ and $\mu_v^-$.  In other words, we have equality of measures $\mu_v^+ = \mu_v^-$ at all places $v$ of $k$.  In particular, letting $v$ be an archimedean place of $k$, we have $\mu_+ = \mu_-$, contradicting Lemma \ref{currents}.
\qed

\begin{remark}
The equidistribution of the postcritically-finite maps in $\Per_1(0)$ follows from Proposition \ref{S equidistribution} in the proof of Theorem \ref{preperiodic}.
\end{remark}

\bigskip
\section{From coincidence to an algebraic relation}
\label{preperiodic proof}

In this section, we complete the proof of Theorem \ref{preperiodic}.  The implications (3) $\implies$ (2) $\implies$ (1) are covered by Proposition \ref{preperiodic easy}.  Throughout this section, we assume condition (1).  We combine the arithmetic equidistribution theorem (Theorem \ref{thm:arithmetic_equidistribution}) with techniques from complex analysis to obtain (3).

\subsection{Preliminary definitions}
Let $G_t$ denote the escape-rate function for $f_t$, as defined in (\ref{escape}), and set 
\begin{equation} \label{eq:G_i}
G_i(t) = \lim_{n\to\infty} \frac{1}{d^n} \log^+ |f_t(a_i(t))| = G_t(a_i(t)).
\end{equation}
Define the bifurcation measure 
\begin{equation}  \label{eq:mu_i}
\mu_i = \frac{1}{2\pi} \Delta G_i
\end{equation}
on the parameter space; by Proposition \ref{active}, the activity of $a_i$ implies that the measure is nonzero.  In fact, we see from the proof of Proposition \ref{active} that the total mass of $\mu_i$ can be computed by the degree growth of the polynomials $f_t^n(a_i(t))$ as $n\to\infty$.  If we pass to a high enough iterate $f^{N_i}_t(a_i(t))$, then 
	$$\deg_t f^{N_i+n}_t(a_i(t)) = m_id^n$$
for some integer $m_i>0$ and all $n\geq 0$.  Then
	$$G_i(t) = \frac{m_i}{d^{N_i}} \log|t| + O(1)$$
as $t\to\infty$; consequently, the measure $\mu_i$ has total mass $m_i/d^{N_i}$.  

For the remainder of the proof, it will be convenient to replace $a_i$ with its iterate $f_t^{N_i}(a_i(t))$.  We may therefore assume that 
\begin{equation}  \label{eq:m_i}
	\deg_t(f_t^n(a_i(t))) = m_i d^n
\end{equation}
for all $n\geq 0$ and 
\begin{equation} \label{mu_i mass}
	\int_{\C} \mu_i \; = \; m_i \, .
\end{equation}

\subsection{(1) $\implies$ ``almost (2)" via arithmetic equidistribution}  \label{arithmetic part}
By assumption, there are infinitely many parameters $t_1,t_2,\ldots\in\C$ such that both $a_1(t_n)$ and $a_2(t_n)$ are preperiodic for $f_{t_n}$.  Following the arguments in \cite{Ghioca:Hsia:Tucker}, Theorem~\ref{thm:arithmetic_equidistribution} implies that the sets
	$$S_1 = \{t:  a_1(t) \mbox{ is preperiodic for } f_t\}$$
and 
	$$S_2 = \{t:  a_2(t) \mbox{ is preperiodic for } f_t\}$$
differ by at most finitely many elements.  If we know that the family $f_t$ and marked points $a_i$ are defined over $\Qbar$, then equidistribution guarantees that $S_1 = S_2$. 

We explain how this follows from \cite{Ghioca:Hsia:Tucker}.  We have already replaced each $a_i$ by a suitably large iterate so that condition (5.1) from \cite{Ghioca:Hsia:Tucker} and the conclusion of their Lemma 5.2 are satisfied for $i=1,2$.  Then for any product formula field $k$ over which $f_t$ and $a_i(t)$ are defined, Corollary 6.11 of \cite{Ghioca:Hsia:Tucker} guarantees that the hypotheses of the equidistribution result Theorem~\ref{thm:arithmetic_equidistribution} are satisfied.  This proves:

\begin{prop}  \label{S equidistribution}
For each $i$, the set $S_i$ is equidistributed with respect to the measure $\mu_{i,v}$ on $\A^1_{\Berk,v}$ for all places $v$ of $k$.  More precisely, given any sequence of subsets $\{E_n\}$ in $S_i$ for which $\# (\Gal(\ksep / k) \cdot E_n) \to \infty$ as $n \to \infty$, the discrete probability measures $[E_n]_v$ converge weakly to $\mu_v$ as $n \to \infty$, for all $v \in \cM_k$.
\end{prop}

\noindent
Consequently, we have $\mu_{1,v} = \mu_{2,v}$ for all places $v$ of $k$.  Here $\mu_{i,v}$ denotes the equilibrium measure on the set $\M_{i,v}$, the closure in $\A^1_{\Berk,v}$ of the set of $t \in \C_v$ for which $a_i(t)$ is bounded under iteration of $f_t$.  It follows that the associated canonical heights $\hat{h}_1$ and $\hat{h}_2$ must be equal.  If $k$ is a number field, the desired equality $S_1 = S_2$ follows, because $S_i = \{t\in \bar{k}: \hat{h}_i(t) = 0\}$ in this case.   The general case follows from Proposition 10.5 of \cite{Ghioca:Hsia:Tucker}.  Note that the hypothesis (i) in Theorem 2.3 of \cite{Ghioca:Hsia:Tucker} is not needed for any of these conclusions.  


\subsection{The boundedness locus $\M$}  \label{phi_M}
Consider the ``generalized Mandelbrot set" associated to $a_i$, defined by
	$$\M_i := \{t \in \C: \mbox{the orbit of } a_i(t) \mbox{ is bounded}\}.$$
As with the usual Mandelbrot set, the boundary of $\M_i$ is the activity locus for $a_i$; that is, the boundary of $M_i$ is the set of parameters $t_0\in\C$ for which $\{t \mapsto f_t^n(a_i(t))\}$ fails to form a normal family on every neighborhood of $t_0$.  The set $S_i$, where $a_i$ is preperiodic, is a subset of $\M_i$.  From Lemma \ref{normal}, the closure of $S_i$ contains the boundary of $\M_i$.  And, exactly as for the usual Mandelbrot set, the Maximum Principle guarantees that the complement of $\M_i$ is connected.  Thus, the conclusion of \S\ref{arithmetic part} (that $S_1$ and $S_2$ differ in at most finitely many elements) guarantees that $\M_1 = \M_2$.  We let $\M$ denote this common set, so
	$$\M := \M_1 = \M_2$$
is the {\em boundedness locus} for $a_1$ and $a_2$.   From Remark \ref{M_a compact}, the set $\M$ is compact.  

Recall that the function $G_i$ defined in (\ref{eq:G_i}) is, up to a multiplicative constant, the Green function for $\M_i$ (see Remark \ref{M_a compact}; cf. \cite[Lemma 6.10]{Ghioca:Hsia:Tucker}).  It follows that 
	$$G_2(t) = \alpha \, G_1(t) \quad\mbox{and} \quad \mu_2 = \alpha \, \mu_1$$
where
\begin{equation} \label{alpha}
	\alpha = \frac{m_2}{m_1} = \frac{\deg a_2(t)}{\deg a_1(t)}
\end{equation}
by equation (\ref{mu_i mass}).  

We will also need the ``uniformizing coordinate" $\phi_M$ associated to the compact set $\M\subset \C$.  This is the uniquely determined univalent function defined in a neighborhood of infinity, with $\phi_M(t) = t + O(1)$ near $\infty$, such that 
	$$\log|\phi_M(t)| = \frac{1}{m_i} \; G_i(t)$$
for $i = 1,2$.  It exists because the periods of the conjugate differential
	$$d^*G_i = -\frac{\del G_i}{\del y} \, dx + \frac{\del G_i}{\del x} \, dy$$
lie in $2\pi m_i \Z$ (for loops around infinity with $t$ large);  see, e.g., \cite[Chapter 4, \S6.1]{Ahlfors}.

\subsection{Analytic relation between $a_1$ and $a_2$}
Let $\phi_t$ denote the uniformizing B\"ottcher coordinate for $f_t$.  That is, for each fixed $t$, $\phi_t$ is defined and univalent in a neighborhood of infinity and is uniquely determined by the conditions that $\phi_t(f_t(z)) = \phi_t(z)^d$ and $\phi_t(z) = z + O(1)$ for all $t$.  
The B\"ottcher coordinate satisfies 
	$$\log|\phi_t(z)| = G_t(z)$$
where it is defined.  See e.g., \cite[\S9]{Milnor:dynamics}.

The following Lemma appears as \cite[Proposition 7.6]{Ghioca:Hsia:Tucker}, but we include a proof for completeness.  The arguments are similar to our proof of \cite[Lemma 3.2]{BD:preperiodic}.

\begin{lemma} \label{domain}
For each $i=1,2$, there exists an integer $n_i$ so that the iterate $f_t^{n_i}(a_i(t))$ lies in the domain of $\phi_t$ for all sufficiently large $t$.
\end{lemma}  

\proof
For each $t$, let $M(f_t)$ denote the maximal critical escape rate, so
	$$M(f_t) = \max\{G_t(c):  f_t'(c) = 0\}.$$
The natural domain of $\phi_t$ is 
	$$\{z\in\C:  G_t(z) > M(f_t)\}.$$
The polynomial growth of the coefficients of $f_t$ implies that $M(f_t)$ grows logarithmically in $t$.  Indeed, by passing to a finite cover of the punctured disk $\{|t| > R\}$ for some $R\gg0$, we may assume that the critical points of $f_t$ are holomorphic functions of $t$.  Applying \cite[Proposition 10.4]{DM:trees}, which uses standard distortion estimates for univalent functions, we conclude that 
	$$M(f_t) = e \log|t| + O(1)$$
as $t\to\infty$ for some $e>0$.  

From Proposition \ref{active}, for each $i$ we know that 
	$$G_t(a_i(t)) = q_i \log|t| + O(1)$$
for some $q_i>0$ as $t\to\infty$.  Choosing $n_i$ so that 
	$$q_i d^{n_i} > e,$$
we conclude that $f_t^{n_i}(a(t))$ lies in the domain of $\phi_t$ for all $t$ sufficiently large. 
\qed

\bigskip
For the rest of this section, we will replace $a_i(t)$ with its iterate $f_t^{n_i}(a_i(t))$ from Lemma \ref{domain}.  

Write each polynomial as 
	$$a_i(t) = \zeta_i t^{m_i} + o(t^{m_i})$$
for some nonzero $\zeta_i\in \C$ and $t$ near infinity.  Define 
	$$\Phi_i (t) = \phi_t(a_i(t))$$
so that 
	$$\Phi_i(t) = \zeta_i t^{m_i} + o(t^{m_i})$$
for $t$ near infinity.  Set 
	$$L = \mathrm{lcm} \{m_1, m_2\}$$
and write 
	$$L = k_1m_1 = k_2 m_2.$$

Now $\Phi_1, \Phi_2$ are analytic maps near infinity, and each satisfies
	$$|\Phi_i(t)| = \exp(G_t(a_i(t)))$$
so that 
	$$\Phi_i(t)  = \zeta_i \, \phi_M(t)^{m_i}$$
for each $i$, where $\phi_M$ is the uniformizing coordinate for $\M$, defined in \S\ref{phi_M}; this is because an analytic map is uniquely determined by its absolute value, up to a rotation.  Therefore, 
	$$\phi_t(a_2(t))^{k_2} = \frac{\zeta_2^{k_2}}{\zeta_1^{k_1}} \, \phi_t(a_1(t))^{k_1}.$$
Set $\zeta = \zeta_2^{k_2} / \zeta_1^{k_1}$.  Then for every $n$, we have 
\begin{equation} \label{k12}
	\phi_t(f_t^n(a_2(t)))^{k_2} = \phi_t(a_2(t))^{k_2d^n} = ( \zeta \phi_t(a_1(t))^{k_1})^{d^n} =  \zeta^{d^n} \phi_t(f_t^n(a_1(t)))^{k_1}.
\end{equation}	
We will refer to (\ref{k12}) as the {\em analytic relation} between the orbits of $a_1(t)$ and $a_2(t)$.  

\begin{lemma} \label{modulus 1}
The $\zeta$ of the analytic relation (\ref{k12}) satisfies $|\zeta| = 1$.
\end{lemma}

\proof
Recall that the constant $\alpha$ from (\ref{alpha}) is given by $\alpha = m_2/m_1$ and the integers $k_1$ and $k_2$ were chosen so that $k_1m_1 = k_2m_2$.  Consequently, 
	$$\log|\Phi_2(t)| = G_t(a_2(t)) = \alpha \, G_t(a_1(t)) = \frac{k_1}{k_2} G_t(a_1(t)) = \frac{k_1}{k_2} \log|\Phi_2(t)|$$
and so 
 	$$|\zeta_2|^{k_2} |\phi_M(t)|^{m_2k_2} = |\Phi_2(t)|^{k_2} = |\Phi_1(t)|^{k_1} =|\zeta_1|^{k_1} |\phi_M(t)|^{m_1k_1}.$$
We see that
	$$|\zeta| = |\zeta_2|^{k_2} / |\zeta_1|^{k_1} = 1,$$
so that 
	$$|\phi_t(a_2(t))|^{k_2} = |\phi_t(a_1(t))|^{k_1}$$
from (\ref{k12}).  
\qed

\subsection{Properties of the B\"ottcher coordinate}
Our next goal will be to promote the analytic relation (\ref{k12}) to an algebraic relation between the orbits of $a_1$ and $a_2$.  To achieve this, we need some estimates on $\phi_t$.  Write the expansion of $\phi_t$ for $z$ near $\infty$ as
	$$\phi_t(z) = z + \sum_{s=1}^\infty g_s(t) z^{-s}.$$
The constant term is 0 because all $f_t$ are centered.  Note that $\phi_t(z)$ is analytic in both $t$ and $z$, where defined.  

\begin{lemma}  \label{polynomial}
The coefficient $g_s(t)$ is polynomial in $t$ for all $s$.  
\end{lemma}

\proof
Recall that  
	$$\phi_t(f^n_t(z)) = (\phi_t(z))^{d^n},$$
for any $n$.  Expand both sides as series in $z$, so we have
	$$f^n_t(z) + O(z^{-d^n}) =  $$  $$z^{d^n} + c g_1(t) z^{d^n-2} + c g_2(t) z^{d^n-3} + (c g_3(t) + c' g_1(t)^2) z^{d^n - 4} + (cg_4(t) + c'' g_1(t) g_2(t)) z^{d^n - 5} + \cdots $$
for nonzero constants $c, c', c'', \ldots$ depending only on $d$ and $n$.  

As the coefficients of the principal part of the left-hand side are polynomials, an induction argument allows us to conclude that the $g_s(t)$ is polynomial for every $s$.  
\qed

\medskip
Let $m = \min\{m_1, m_2\}$, where $m_i$ is the degree in $t$ of $a_i(t)$.

\begin{lemma}  \label{degree bound}
The degree of $g_s(t)$ in $t$ is no greater than $m(s+1)$.  
\end{lemma}

\proof
For fixed $t$, choose $R = R_t$ minimal such that $\{z: |z|>R\}$ lies in the domain of the univalent function $\phi_t$.  Then
	$$\psi(z) = R/\phi_t(R/z)$$
defines a univalent function on the unit disk, with $\psi(0) = 0$ and $\psi'(0) = 1$.  Expand $\psi$ in a power series around 0, so 
	$$\psi(z) = z + \sum_{n=2}^\infty b_n z^n.$$
Littlewood's Theorem implies that $|b_n| \leq e \, n$ for all $n$ and any such $\psi$;  see \cite[\S2.4]{Duren:univalent}.  In fact, by the Bieberbach Conjecture (Theorem of de Branges), we know that $|b_n| \leq n$, but this is not necessary for us.  

In our case, the first few terms in the expansion of $\psi$ are 
	$$\psi(z) = z - g_1 R^{-2} z^3 + g_2 R^{-3} z^4 + (g_1^2 - g_3)R^{-4} z^5 + (2g_1g_2 - g_4) R^{-5} z^6 + \cdots$$
An induction argument implies that 
	$$|g_s| \leq C_s R^{s+1}$$
for some constant depending on $s$, where $g_s$ and $R$ both depend on $t$. 

Now, recall that $a_i(t)$ lies in the domain of $\phi_t$ for all $t$ sufficiently large, by Lemma \ref{domain} (and the comment following the proof).  From the traditional distortion arguments applied to $\phi_t$ (applying the estimate $|b_2| \leq 2$; see \cite[Lemma 3.2]{BD:preperiodic} or \cite[\S3]{Branner:Hubbard:1}), the region 
	$$\{z: G_t(z) > d\, M(f_t)\}$$
lies outside the disk of radius $R_t$ for all $t$ large.  And so we may assume that $a_i(t)$ lies outside the disk of radius $R_t$ for all $t$ large.  That is, we have $R_t \leq |a_i(t)|$ for $|t|\gg0$, and we conclude that
	$$|g_s(t)| \leq C_s |a_i(t)|^{s+1}$$	
for $i=1,2$.  Finally, then, the degree of $g_s$ must be no greater than the degree of $a_i(t)^{s+1}$.  
\qed

\subsection{Polynomial relation between $a_1(t)$ and $a_2(t)$}  \label{relation section}
Expand each power of the B\"ottcher coordinate $\phi_t(z)$ in Laurent series near infinity as
	$$(\phi_t(z))^k = P_t^k(z) + \sum_{s=1}^\infty b^k_s(t) z^{-s}.$$
By Lemma \ref{polynomial}, the expression $P_t^k(z)$ is polynomial in both $t$ and $z$; in $z$ it is monic and centered of degree $k$.  By Lemma \ref{degree bound}, we may conclude that 
\begin{equation} \label{b degree}
	\deg_t b^k_s(t) \leq m(s+k).
\end{equation}
Indeed, $b_s^k$ is a sum of products $\prod_{i=1}^l g_{s_i}$ for some $l\in \{1, \ldots, k\}$ where $\sum s_i = k-l+s$, and each product has degree at most $\sum_{i=1}^l m(s_i+1) = ml + m\sum s_i = ml + mk - ml + ms = m(s+k)$.

Setting $z = f_t^n(a_i(t))$, we have 
	$$(\phi_t(f_t^n(a_i(t))))^k = P_t^k(f_t^n(a_i(t))) + \sum_{s=1}^\infty b^k_s(t) (f_t^n(a_i(t)))^{-s}.$$
By (\ref{b degree}), the infinite-sum term is $O(t^{-m(d^n - k -1)})$.  Recall from equation (\ref{k12}) we have
	$$\phi_t(f_t^n(a_2(t)))^{k_2} =  \zeta^{d^n} \phi_t(f_t^n(a_1(t)))^{k_1}$$
for all $n\geq 0$.   Expanding both sides in $t$ implies that the polynomial parts of both sides must be equal for any $n$.  Thus, for all $n\gg0$, we have 
\begin{equation} \label{poly relation}
	P^{k_2}_t(f_t^n(a_2(t))) = \zeta^{d^n} P^{k_1}_t(f_t^n(a_1(t))).
\end{equation}
It will be convenient to replace $a_1$ and $a_2$ with higher iterates so that equation (\ref{poly relation}) holds for {\em all} $n$.  

We would like to know that the polynomial relation (\ref{poly relation}) between $f_t^n(a_2(t))$ and $f_t^n(a_1(t))$ is independent of $n$, or at least that the constants $\zeta^{d^n}$ cycle through only finitely many values.  We thank Dragos Ghioca for suggesting the strategy for this proof of Lemma \ref{general zeta}.

\begin{lemma}  \label{general zeta}
The ratio $\zeta = \zeta_2^{k_2}/\zeta_1^{k_1}$ is a root of unity.  
\end{lemma}

\proof  
We combine the analytic relation (\ref{k12}) and the polynomial relation (\ref{poly relation}) to obtain
\begin{equation} \label{sum}
	\sum_{s=1}^\infty b_s^{k_2}(t) (f_t^n(a_2(t)))^{-s} = \zeta^{d^n} \sum_{s=1}^\infty b_s^{k_1}(t) (f_t^n(a_1(t)))^{-s}
\end{equation}
for all $n$.  Let $s_1$ be the smallest $s$ for which $b_s^{k_1}$ is nonzero and $s_2$ the smallest $s$ for which $b_s^{k_2}$ is nonzero.  Let $C(s,k)$ denote the degree of the polynomial $b_s^k(t)$; recall from (\ref{b degree}) that $C(s,k) \leq s(m+k)$.  Expanding both sides of (\ref{sum}) in $t$, the leading term on the left-hand-side is  
	$$c_2 \zeta_2^{-s_2d^n} t^{C(s_2,k_2) - s_2m_2d^n}$$
for some constant $c_2\in\C^*$, while the leading term on the right-hand-side is 
	$$c_1 \zeta^{d^n} \zeta_1^{-s_1d^n} t^{C(s_1, k_1) - s_1m_1d^n}$$
for some constant $c_1\in \C^*$.  As we have equality in (\ref{sum}) for all $n$, it follows that $s_2m_2 = s_1m_1$.   As $L = k_1m_1 = k_2m_2$ is the least common multiple of $m_1$ and $m_2$, we may write $s_1 = \ell k_1$ and $s_2 = \ell k_2$ for some positive integer $\ell$. Furthermore, the coefficients of the leading terms must coincide, so 
	$$\frac{c_2}{c_1} = \zeta^{d^n} \left(\frac{\zeta_2^{k_2}}{\zeta_1^{k_1}} \right)^{\ell d^n} = \zeta^{d^n + \ell d^n}$$
for all $n$.  Therefore, $\zeta$ is a root of unity.
\qed 

\begin{remark}  The proof of Lemma \ref{general zeta} is elementary but somewhat unenlightening.  
When the points $a_i(t)$ are critical (i.e., in the setting of Theorem \ref{main}), one can give a more conceptual proof that $\zeta$ is a root of unity, as follows.
From Lemma \ref{modulus 1}, we know that $|\zeta|= 1$.  From (\ref{k12}), the argument of $\zeta$ is equal to the difference in argument between $\phi_t(a_2(t))^{k_2}$ and $\phi_t(a_1(t))^{k_1}$, independent of $t$.  
We are assuming that there are infinitely many parameters $t$ such that $f_t$ is PCF, and all periodic cycles for a PCF map must be superattracting or repelling \cite[Corollary 14.5]{Milnor:dynamics}.  From \S\ref{arithmetic part} and Lemma \ref{normal}, there are infinitely many $t\in\del \M$ such that both $a_1(t)$ and $a_2(t)$ are preperiodic to repelling cycles.
Such a parameter $t_0$ will be a landing point of a ``rational external ray" for $\phi_\M$ (see e.g. \cite[Chapter 18]{Milnor:dynamics}).
In other words, the points $a_1(t_0)$ and $a_2(t_0)$ will be landing points for rational external rays in the Julia set of $f_{t_0}$.  
It follows that the difference in argument between $\phi_t(a_1(t))$ and $\phi_t(a_2(t))$ is rational, and therefore that $\zeta$ is a root of unity.
\end{remark}

\bigskip
Lemma \ref{general zeta} implies that the sequence $\{\zeta^{d^n}: n\geq 0\}$ will eventually cycle.  Replacing $a_1$ and $a_2$ with iterates will allow us to assume that $\zeta$ itself is periodic for $z^d$.  That is, we may assume there exists a positive integer $k$ so that 
	$$\zeta ^{d^k} = \zeta.$$
Equation (\ref{poly relation}) can be formulated as
\begin{equation} \label{invariant}
	P^{k_2}_t(f_t^{kn}(a_2(t))) = \zeta \, P^{k_1}_t(f_t^{kn}(a_1(t)))
\end{equation}
for all $n$.

\subsection{Simplifying the algebraic relation (\ref{invariant}) and concluding the proof}  \label{conclusion}
Define polynomials 
	$$A_t(z) := P^{k_1}_t(z) \qquad B_t(z) := \zeta P^{k_2}_t(z).$$
Then (\ref{invariant}) implies that the algebraic curve (or a subset of the irreducible components of this curve)
	$$\{(x,y):  A_t(x) = B_t(y)\} \subset \P^1\times\P^1$$
is invariant for the map 
	$$(f_t^k, f_t^k) : \P^1\times\P^1 \to \P^1 \times \P^1$$
for every $t$.

If the polynomial $A_t(x) - B_t(y)$ is reducible for all $t$, let $Q_t(x,y)$ denote a factor such that $Q_t(a_1(t), a_2(t)) = 0$ and $Q_t$ is irreducible for general $t$.   There are only finitely many irreducible components, so by passing to higher iterates (of the $a_i$ and of the $f^k$ preserving $\zeta$), we may assume that the curve 
	$$C_t = \{Q_t(x,y) = 0\}$$ 
is invariant for $(f_t^k, f_t^k)$ for all $t$, and $C_t$ is irreducible for general $t$.  

We now appeal to the classification of $(f,f)$-invariant curves in $\P^1\times\P^1$ for polynomials $f$.  It was treated in great generality by Medvedev and Scanlon in \cite{Medvedev:Scanlon}, applying Ritt's study of polynomial decompositions from \cite{Ritt:decompositions}.  As the family $f_t$ is nontrivial, the polynomial $f_t^k$ cannot be conjugate to $z^{d^k}$ or a Chebyshev polynomial for all $t$.  Therefore, the curve $C_t$ must be a graph, of the form $\{y = h_t(x)\}$ or $\{x = h_t(y)\}$, for a polynomial $h_t$ that commutes with $f_t^k$ \cite[Theorem 6.24]{Medvedev:Scanlon}.  

In other words, there exists a polynomial $h\in \C[t,z]$ such that $h_t\circ f_t^k = f_t^k\circ h_t$ for all $t$ and so that either $a_2(t) = h_t(a_1(t))$ or $a_1(t) = h_t(a_2(t))$ for all $t$.  (Recall that we have repeatedly replaced the original $a_i$ with an iterate $f_t^{n_i}(a_i(t))$.)  If the conclusion is that $a_1(t) = h_t(a_2(t))$, then the proof of condition (3) is complete.  Suppose instead that $a_2(t) = h_t(a_1(t))$.  If $\deg_z h = 1$, then we may replace $h_t$ with $h_t^{-1}$ to achieve the conclusion of condition (3).  If $\deg_z h_t > 1$, then from Ritt's work we know that $h_t$ must share an iterate with $f_t$; say $h_t^q = f_t^r$ \cite{Ritt:decompositions, Ritt:permutable}.  Then $h_t^{q-1}(a_1(t)) = f_t^r(a_2(t))$, so we again achieve the conclusion of condition (3), taking the new $h$ to be $h_t^{q-1}(z)$.   This concludes the proof.

\bigskip
\section{Proof of Theorem \ref{main}}
In this final section, we provide the proof of Theorem \ref{main}.  In most respects, Theorem \ref{main} is a special case of Theorem \ref{preperiodic}.  

\medskip{\bf (1) $\implies$ (2)}.
Let $a_1(t)$ and $a_2(t)$ denote any pair of active critical points of $f_t$.  At each postcritically-finite polynomial $f_t$, both $a_1(t)$ and $a_2(t)$ have finite forward orbit.  From Theorem \ref{preperiodic}, condition (1) implies that the sets $S_1$ and $S_2$ coincide.  In addition, as observed in \S\ref{phi_M}, the sets $\M_1$ and $\M_2$ must coincide, and therefore so do their harmonic measures (relative to $\infty$).  From Remark \ref{M_a compact}, the harmonic measure on $\M_i$ is exactly the bifurcation measure for the critical point $a_i$, normalized to have total mass 1. 

\medskip{\bf (2) $\implies$ (3)}.
For each active critical point $c_i$, the support of the bifurcation measure $\mu_i$ is equal to the (outer) boundary of the set $\M_i$.  Each $\M_i$ is full (meaning that its complement is connected):  indeed, on a bounded component of $\C\setminus \M_i$, the Maximum Principle guarantees that the magnitude of $f_t^n(c_i(t))$ never exceeds its maximum value on $\M_i$.  Therefore the measure $\mu_i$ determines $\M_i$ as a set.  And so $\M_i$ does not depend on the choice of active critical point.  In particular, all critical points have bounded forward orbit for $f_t$ if and only if $t\in \M_i$ for some active critical point $i$.  Therefore, $\M_i$ is the connectedness locus for $f_t$.  

\medskip{\bf (3) $\implies$ (4)}.  This implication is exactly as in the proof of Theorem \ref{preperiodic}.  Specifically, the arguments of \S\ref{phi_M}--\ref{conclusion} start with the assumption that the sets $\M_i$ coincide and conclude with the desired algebraic relation (4).

\medskip{\bf (4) $\implies$ (1)}.  If a critical point is not active, then Proposition \ref{active} shows that it is preperiodic for all parameters $t\in\C$.  If there is only one active critical point, then Lemma \ref{normal} implies that it has finite orbit for infinitely many $t$, and therefore $f_t$ is postcritically finite for infinitely many $t$.  If there are at least two active critical points, then Theorem \ref{preperiodic} guarantees that they are simultaneously preperiodic at infinitely many parameters $t$.  Again we conclude that $f_t$ is postcritically finite for infinitely many $t$.

 \bigskip\bigskip

\begin{thebibliography}{GHT}

\bibitem[Ah]{Ahlfors}
L~V. Ahlfors.
\newblock {\em Complex analysis}.
\newblock McGraw-Hill Book Co., New York, third edition, 1978.
\newblock An introduction to the theory of analytic functions of one complex
  variable, International Series in Pure and Applied Mathematics.

\bibitem[An]{Andre:finitude}
Y.~Andr{\'e}.
\newblock {Finitude des couples d'invariants modulaires singuliers sur une
  courbe alg\'ebrique plane non modulaire}.
\newblock {\em J. Reine Angew. Math.} {\bf 505}(1998), 203--208.

\bibitem[BD]{BD:preperiodic}
M.~Baker and L.~DeMarco.
\newblock {Preperiodic points and unlikely intersections}.
\newblock {\em Duke Math. J.} {\bf 159}(2011), 1--29.

\bibitem[BR]{BRBook}
M.~Baker and R.~Rumely.
\newblock {\em Potential theory and dynamics on the {B}erkovich projective
  line}, volume 159 of {\em Mathematical Surveys and Monographs}.
\newblock American Mathematical Society, Providence, RI, 2010.

\bibitem[Be1]{Beardon:symmetries}
A.~F. Beardon.
\newblock {Symmetries of {J}ulia sets}.
\newblock {\em Bull. London Math. Soc.} {\bf 22}(1990), 576--582.

\bibitem[Be2]{Benedetto:polynomial}
R.~L. Benedetto.
\newblock {Heights and preperiodic points of polynomials over function fields}.
\newblock {\em Int. Math. Res. Not.} {\bf 62}(2005), 3855--3866.

\bibitem[BH1]{Branner:Hubbard:1}
B.~Branner and J.~H. Hubbard.
\newblock {The iteration of cubic polynomials. {I}. {T}he global topology of
  parameter space}.
\newblock {\em Acta Math.} {\bf 160}(1988), 143--206.

\bibitem[BE]{Buff:Epstein:PCF}
X.~Buff and A.~Epstein.
\newblock {Bifurcation measure and postcritically finite rational maps}.
\newblock In {\em Complex dynamics}, pages 491--512. A K Peters, Wellesley, MA,
  2009.

\bibitem[CL]{ChambertLoir}
A.~Chambert-Loir.
\newblock {Mesures et \'equidistribution sur les espaces de {B}erkovich}.
\newblock {\em J. Reine Angew. Math.} {\bf 595}(2006), 215--235.

\bibitem[De1]{D:current}
L.~DeMarco.
\newblock {Dynamics of rational maps: a current on the bifurcation locus}.
\newblock {\em Math. Res. Lett.} {\bf 8}(2001), 57--66.

\bibitem[De2]{D:lyap}
L.~DeMarco.
\newblock {Dynamics of rational maps: {L}yapunov exponents, bifurcations, and
  capacity}.
\newblock {\em Math. Ann.} {\bf 326}(2003), 43--73.

\bibitem[DM]{DM:trees}
L.~DeMarco and C.~McMullen.
\newblock {Trees and the dynamics of polynomials}.
\newblock {\em Ann. Sci. {\'E}cole Norm. Sup.} {\bf 41}(2008), 337--383.

\bibitem[DH]{Douady:Hubbard:Thurston}
A.~Douady and J.~H.~.Hubbard.
\newblock{A proof of Thurston's topological characterization of rational functions}.
\newblock {\em  Acta Math.} {\bf 171}(1993), no. 2,  263--297.

\bibitem[DF]{Dujardin:Favre:critical}
R.~Dujardin and C.~Favre.
\newblock {Distribution of rational maps with a preperiodic critical point}.
\newblock {\em Amer. J. Math.} {\bf 130}(2008), 979--1032.

\bibitem[Du]{Duren:univalent}
P.~L. Duren.
\newblock {\em Univalent functions}, volume 259 of {\em Grundlehren der
  Mathematischen Wissenschaften [Fundamental Principles of Mathematical
  Sciences]}.
\newblock Springer-Verlag, New York, 1983.

\bibitem[FRL]{FRL:equidistribution}
C.~Favre and J.~Rivera-Letelier.
\newblock {\'{E}quidistribution quantitative des points de petite hauteur sur
  la droite projective}.
\newblock {\em Math. Ann.} {\bf 335}(2006), 311--361.

\bibitem[GHT1]{Ghioca:Hsia:Tucker}
D.~Ghioca, L.-C. Hsia, and T.~Tucker.
\newblock {Preperiodic points for families of polynomials}.
\newblock {\em {\em To appear}, Algebra \& Number Theory}.

\bibitem[GHT2]{GHT:2012}
D.~Ghioca, L.-C. Hsia, and T.~Tucker.
\newblock {Preperiodic points for families of rational maps}.
\newblock {Submitted for publication, 2012}.

\bibitem[GTZ]{Ghioca:Tucker:Zhang}
D.~Ghioca, T.~Tucker, and S.~Zhang.
\newblock {Towards a dynamical {M}anin-{M}umford conjecture}.
\newblock {\em Int. Math. Res. Not.} {\bf 22}(2011), 5109--5122.

\bibitem[Hi]{Hindry}
M.~Hindry.
\newblock {Autour d'une conjecture de {S}erge {L}ang}.
\newblock {\em Invent. Math.} {\bf 94}(1988), 575--603.

\bibitem[La]{Laurent}
M.~Laurent.
\newblock {\'{E}quations diophantiennes exponentielles}.
\newblock {\em Invent. Math.} {\bf 78}(1984), 299--327.

\bibitem[MSS]{Mane:Sad:Sullivan}
R.~Ma\~{n}{\'e}, P.~Sad, and D.~Sullivan.
\newblock {On the dynamics of rational maps}.
\newblock {\em Ann. Sci. Ec. Norm. Sup.} {\bf 16}(1983), 193--217.

\bibitem[Mc1]{McMullen:families}
C.~T. McMullen.
\newblock {Families of rational maps and iterative root-finding algorithms}.
\newblock {\em Ann. of Math. (2)}, {\bf 125}(1987), 467--493.

\bibitem[Mc2]{McMullen:universal}
C.~T. McMullen.
\newblock {The {M}andelbrot set is universal}.
\newblock In {\em The {M}andelbrot set, theme and variations}, volume 274 of
  {\em London Math. Soc. Lecture Note Ser.}, pages 1--17. Cambridge Univ.
  Press, Cambridge, 2000.

\bibitem[MS]{Medvedev:Scanlon}
A.~Medvedev and T.~Scanlon.
\newblock {Invariant varieties for polynomial dynamical systems}.
\newblock {To appear, \em Ann. of Math.}.

\bibitem[Mi1]{Milnor:cubicpoly}
J.~Milnor.
\newblock {Remarks on iterated cubic maps}.
\newblock {\em Experiment. Math.} {\bf 1}(1992), 5--24.

\bibitem[Mi2]{Milnor:dynamics}
J.~Milnor.
\newblock {\em Dynamics in one complex variable}, volume 160 of {\em Annals of
  Mathematics Studies}.
\newblock Princeton University Press, Princeton, NJ, {T}hird edition, 2006.

\bibitem[Pi]{Pila:AO}
J.~Pila.
\newblock {O-minimality and the {A}ndr\'e-{O}ort conjecture for {$\mathbb C^n$}}.
\newblock {\em Ann. of Math. (2)} {\bf 173}(2011), 1779--1840.

\bibitem[Ra1]{Raynaud:1}
M.~Raynaud.
\newblock {Courbes sur une vari\'et\'e ab\'elienne et points de torsion}.
\newblock {\em Invent. Math.} {\bf 71}(1983), 207--233.

\bibitem[Ra2]{Raynaud:2}
M.~Raynaud.
\newblock {Sous-vari\'et\'es d'une vari\'et\'e ab\'elienne et points de
  torsion}.
\newblock In {\em Arithmetic and geometry, {V}ol. {I}}, volume~35 of {\em
  Progr. Math.}, pages 327--352. Birkh\"auser Boston, Boston, MA, 1983.

\bibitem[Ri1]{Ritt:decompositions}
J.~F. Ritt.
\newblock {Prime and composite polynomials}.
\newblock {\em Trans. Amer. Math. Soc.} {\bf 23}(1922), 51--66.

\bibitem[Ri2]{Ritt:permutable}
J.~F. Ritt.
\newblock {Permutable rational functions}.
\newblock {\em Trans. Amer. Math. Soc.} {\bf 25}(1923), 399--448.

\bibitem[Si]{Silverman:moduli}
J.~H. Silverman.
\newblock {\em Moduli spaces and arithmetic dynamics}, volume~30 of {\em CRM
  Monograph Series}.
\newblock American Mathematical Society, Providence, RI, 2012.




\end{thebibliography}
\def\cprime{$'$}

 \bigskip

\end{document}